\documentclass[2p]{article}
\usepackage{amssymb,amsmath,amsthm}
\usepackage{indentfirst}
\usepackage{exscale}
\usepackage{relsize}
\usepackage[numbers,sort&compress]{natbib}

\usepackage{geometry}
\usepackage{color}
\newcommand{\R}{\mathbb{R}}
\newtheorem{theorem}{Theorem}[section]

\theoremstyle{definition}

\newtheorem{example}[theorem]{Example}

\newtheorem{lemma}[theorem]{Lemma}

\allowdisplaybreaks
\theoremstyle{remark}
\newtheorem{remark}[theorem]{Remark}
\numberwithin{equation}{section}
 \UseRawInputEncoding

\begin{document}
\title{\Large\bf{Nontrivial solutions for a generalized poly-Laplacian system on finite graphs}
 }
\date{}

\author {\ Wanting Qi$^{1}$, Xingyong Zhang$^{1,2}$\footnote{Corresponding author, E-mail address: zhangxingyong1@163.com}\\
      {\footnotesize $^{1}$Faculty of Science, Kunming University of Science and Technology, Kunming, Yunnan, 650500, P.R. China.}\\
      {\footnotesize $^{3}$Research Center for Mathematics and Interdisciplinary Sciences, Kunming University of Science and Technology,}\\
 {\footnotesize Kunming, Yunnan, 650500, P.R. China.}\\
 }

 \date{}
 \maketitle

 \begin{center}
 \begin{minipage}{15cm}
 \par
 \small  {\bf Abstract:} We investigate the existence and multiplicity of solutions for a class of generalized  coupled system involving poly-Laplacian and a parameter $\lambda$ on finite graphs. By using mountain pass lemma together with cut-off technique, we obtain that system has at least a nontrivial weak solution $(u_{\lambda},v_{\lambda})$  for every large  parameter $\lambda$ when the nonlinear term $F(x,u,v)$ satisfies superlinear growth conditions only in a neighborhood of origin point $(0,0)$. We also obtain a concrete form for the lower bound of parameter $\lambda$ and the trend of $(u_{\lambda},v_{\lambda})$  with the change of parameter $\lambda$.
 Moreover, by using a revised Clark's theorem together with cut-off technique, we obtain that system has a sequence of solutions tending to 0 for every $\lambda>0$ when the nonlinear term $F(x,u,v)$ satisfies sublinear growth conditions only in a neighborhood of origin point $(0,0)$.
 \end{minipage}
 \end{center}
  \allowdisplaybreaks
 \vskip2mm
 {\section{Introduction and main results }}
The existence and multiplicity of solutions for the quasilinear elliptic system with the following form
\begin{equation}\label{d1}
 \begin{cases}
 -\Delta_{p}u+h_{1}(x)|u|^{p-2}u=\lambda F_{u}(x,u,v),  \;\;\; x \in \Omega,\\
 -\Delta_{q}v+h_{2}(x)|v|^{q-2}v=\lambda F_{v}(x,u,v),  \;\;\; x \in \Omega
   \end{cases}
 \end{equation}
have been studied extensively over the past several decades (for example, see \cite{Boccardo-Figueiredo-2002,Djellit2003,Djellit2004,Djellit2007} and references therein), where $\Omega\subset\R^{N}$ is a domain, $\Delta_{p}$ is the so-called $p$-Laplacian operator, i.e., $\Delta_{p}=\mbox{div}(|\nabla u|^{p-2}\nabla u)$, $1<p,q<N$, $\lambda>0$ is a parameter, $h_{i}(x)$, $i=1,2$ and $F$ are real valued functions satisfying some suitable assumptions. Specially, in \cite{Boccardo-Figueiredo-2002}, Boccardo investigated the system (\ref{d1}) with $\lambda =1$ and $h_{i}(x)=0$, $i=1,2$, where $\Omega$ is a bounded domain in $\R^{N}$. They assumed that $F:\bar{\Omega}\times\R\times\R\rightarrow \R$ is $C^{1}$ and satisfies
$$  |F(x,u,v)|\leq C(1+|u|^{r}+|v|^{s}),$$
where $C>0$ is a constant, $r\leq p^{*}$, $s\leq q^{*}$, $p^{*}=\frac{pN}{N-p}$ and $q^{*}=\frac{qN}{N-q}$. Moreover, there exist constants $R>0$, $\theta_{p}$ and $\theta_{q}$ with $\frac{1}{p^{*}}<\theta_{p}<\frac{1}{p}$ and $\frac{1}{q^{*}}<\theta_{q}<\frac{1}{q}$ such that
$$0<F(x,u,v)\leq \theta_{p}uF_{u}(x,u,v)+ \theta_{q}vF_{v}(x,u,v)$$
for all $x\in \bar{\Omega}$ and $|u|,|v|\geq R $, which is usually called as Ambrosseti-Rabinowitz condition ((AR) for short) and
$$  |F_{u}(x,u,v)|\leq C(1+|u|^{p^{*}-1}+|v|^{\frac{q^{*}(p^{*}-1)}{p^{*}}}),$$
$$  |F_{v}(x,u,v)|\leq C(1+|v|^{q^{*}-1}+|u|^{\frac{p^{*}(q^{*}-1)}{q^{*}}})$$
and some other reasonable conditions. Under the above conditions, they established the existence results of a nontrivial solution for system (\ref{d1}) by using mountain pass theorem.
\par
Corresponding to system (\ref{d1}), the following well-known $p$-Laplacian equation
\begin{equation}\label{d2}
-\Delta_{p} u +h(x)|u|^{p-2}u
     =\lambda f(x,u),  \;\;\;  x \in \Omega
   \end{equation}
has also attracted great attention as a result of its important  applications to many fields, such as non-Newtonian fluids, reaction-diffusion problems, elastic mechanics, flow through porous media, thermal radiation and so on. For details, one can see \cite{D-J-I-1985} and references therein. Specially, one of the applications of equation (\ref{d2}) is to describe the steady-state flow behavior of fluid with flow density $u$ in the presence of interference, where $h(x)|u|^{p-2}u$ and $\lambda f(x,u)$ represent two interference sources, respectively. After the interference source $h(x)|u|^{p-2}u$ disappears, (\ref{d2}) reduces to the following $p$-Laplacian elliptic equation
\begin{equation}\label{d3}
-\Delta_{p} u
     =\lambda f(x,u),  \;\;\;  x \in \Omega,
   \end{equation}
which associated with non-Newtonian seepage phenomena \cite{D-P-2007}. At present, there are many achievements in the study of existence of solutions to (\ref{d2}) and
(\ref{d3}) (for example, see \cite{Li-Zhou-2002,Bonanno-Bisci2010,Ramaswamy-Shivaji2004,Costa-Wang2005} and references therein). In most of these papers, nonlinear term is required to have growth both near $0$ and  the infinity  about $u$. In  \cite{Costa-Wang2005}, Costa and Wang make some assumptions on the nonlinear term referred  only to its behavior in a vicinity of $u=0$. By a cut-off technique together with some estimates for solutions as the form $\|u\|_{\infty}\leq c_1\lambda^{-c_2}$ with $c_i>0$, $i=1,2$, they obtained the number of signed and sign-changing solutions for the elliptic problems (\ref{d3}) with $p=2$ and $f(x,u)=f(u)$ under the Dirichlet boundary condition $u(x)=0$, $x\in \partial \Omega$, where $\Omega\subset \R^N$ is bounded domain. The idea in \cite{Costa-Wang2005} has been applied to a lot of different problems, for example,  Schr\"odinger equation(\cite{Huang-Jia2019,Medeiros2008}),  quasilinear elliptic problems involving $p$-Laplacian operator in $\R^{N}$(\cite{Medeiros2007}) and  fractional order Kirchhoff-Type system(\cite{Kang-Liu2020}). Especially, in \cite{Kang-Liu2020}, Kang, Liu and Zhang applied the idea  to the fractional order Kirchhoff-Type system with local superlinear nonlinearity. By using the mountain pass theorem, they established the existence result of a nontrivial weak solution for every given sufficiently large $\lambda $. More importantly, they obtained a concrete value of the lower bound of the parameter $\lambda $ and two estimates of weak solutions.
\par
Recently, the subjects about the partial differential equations on graphs have attracted considerable attentions because of their applications to many fields including data analysis, image processing, neural network, transportation  and so on (see \cite{Chakik2014,Curtis1990,Desquesnes2013,Elmoataz2012,Elmoataz2017,Elmoataz2015} and references therein) and there has been some works, for instance, the existence and non-existence of solutions of heat equations (\cite{Chung2011,Lin2017,Liu2016}), the blow-up and uniqueness properties for solutions of heat equations (\cite{Huang2012,Xin2014}), the existence of solutions of Kazdan-Warner equations (\cite{Grigor2016,Ge2018}). Especially, the existence and multiplicity of solutions for the poly-Laplacian system with the following form
\begin{equation}\label{aa1}
 \begin{cases}
\pounds_{m_{1},p_{1}}u+h_{1}(x)|u|^{p_{1}-2}u=\lambda F_{u}(x,u,v),  \;\;\; x \in V,\\
\pounds_{m_{2},p_{2}}v+h_{2}(x)|v|^{p_{2}-2}v=\lambda F_{v}(x,u,v),  \;\;\; x \in V,
   \end{cases}
 \end{equation}
have also been investigated extensively (for example, see \cite{Yu-Xie-Zhang-Qi,Zhang-Zhang-Xie-Yu-2022,Yu-Zhang-Xie-Zhang-2023,Yang-Zhang-2023} and references therein), where
$G=(V,E)$ is a finite graph, $\pounds_{m_{i},p_{i}}$ are the poly-Laplacian operators  defined in section 2, $m_{i}\geq 1$ and $p_{i}> 1$, $i=1,2$, $h_{i}$ and $F$ are real valued functions satisfying some assumptions.
In \cite{Yu-Xie-Zhang-Qi}, by the variational methods which are based on mountain pass theorem and topological degree theory, the authors obtained that system (\ref{aa1}) has infinitely many solutions.
In \cite{Zhang-Zhang-Xie-Yu-2022}, when $F$ satisfies the super-$(p,q)$-linear growth condition, by applying the mountain pass theorem and the symmetric mountain pass theorem, the authors proved that the system (\ref{aa1}) has at least one nontrivial solutions and $\mbox{dim }W $ nontrivial solutions, where $W$ is the working space of system  (\ref{aa1}), respectively.
In \cite{Yu-Zhang-Xie-Zhang-2023}, when $F$ satisfies the asymptotically-$p$-linear growth conditions at infinity, the authors obtained some sufficient conditions about the existence of a nontrivial solution for system (\ref{aa1}) via the mountain pass theorem.
In \cite{Yang-Zhang-2023},  by using the mountain pass theorem and Ekeland's variational principle, the authors showed that the existence of two nontrivial solutions for system  (\ref{aa1}) involving concave-convex nonlinearity.

\par
In system (\ref{aa1}), if we set $m_{1}=m_{2}$, $p_{1}=p_{2}$, $h_{1}=h_{2}$, $\lambda=1$ and $u=v$,  it reduces to the following elliptic equation:
\begin{equation}\label{d4}
\pounds_{m,p}u+h(x)|u|^{p-2}u
     = f(x,u),
   \end{equation}
there are also many papers studying  problem (\ref{d4}), for example, see \cite{Grigor-Lin2016,Grigor-Lin2017,Zhang-Zhao2018,Han-Shao2021} and references therein.
 In \cite{Grigor-Lin2016}, Grigor'yan, Lin and Yang studied the existence of nontrivial solution for the problem (\ref{d4}) on finite graph by variational method. They assumed that $f(x,t)$ satisfies the following conditions:
\vskip2mm
\noindent
(f1) there exists some $q > p > 1$ and $M > 0$  such that for all $x \in V$ and $|t| \geq M$,
$$
0 < qF(x,t) \leq tf(x,t);
$$
(f2) $\limsup_{t \rightarrow 0} \frac{|f(x,t)|}{|t|^{p-1}}< \lambda_{mp} $, where $\lambda_{mp}$ is the first eigenvalue of $-\Delta_{p}$ defined as
$$\lambda_{mp}  = \inf_{u \not\equiv 0} \frac{\int_{V} (|\nabla^{m} u|^{p} + h|u|^{p}) d\mu}{\int_{V} |u|^{p} d\mu}.$$
They proved that system (\ref{d4}) admits a nontrivial solution by the mountain pass theorem. In \cite{Grigor-Lin2017}, by the same method, they established the existence result of strictly positive solutions for nonlinear equations (\ref{d4}) with $m=1$ and $p=2$. They assumed that $h(x)$ satisfies some reasonable assumptions and nonlinear term $f(x,t)$ satisfies the following conditions:
\vskip2mm
\noindent
(h1) $f(x,0) = 0$, $f(x,t)$ is continuous in $t$ and for any fixed $S > 0$, there exists a constant $M_{S}$ such that $\max_{t \in [0,S]} f(x,t) \leq M_{S}$ for all $x \in V$;\\
(h2) there exists a constant $\mu > 2$ such that for all $x \in V$ and $t\in(0,+\infty)$,
$$
0 <\mu F(x,t)= \mu\int_{0} ^{t} f(x,s) ds \leq tf(x,t) ;
$$
(h3) $\limsup_{t \rightarrow 0^{+}} \frac{2F(x,t)}{|t|^{2}} < \lambda_{1} $, where $\lambda_{1}$ is the first eigenvalue of $-\Delta$ defined as
$$\lambda_{1} = \inf_{\int_{V} u^{2} d\mu = 1} \int_{V} (|\nabla u|^{2} + h|u|^{2}) d\mu.$$
Moreover, in \cite{Zhang-Zhao2018}, by using the Nehari method, on a locally finite graph $G=(V,E)$, Zhang and Zhao investigated the convergence of ground state solutions for the equation (\ref{d4}) with $m=1$, $p=2$, $h(x)=\varsigma a(x)+1$ and $f(x,u)=|u|^{\alpha-2}u$, where $\varsigma$ is a parameter, $\alpha\geq 2$ and $a(x):V\rightarrow [0,+\infty)$. Subsequently,  Han and Shao \cite{Han-Shao2021} further studied the nonlinear $p$-Laplacian equations (\ref{d4}) with $m=1$, $p\geq2$ and $h(x)=\varsigma a(x)+1$ on a locally finite graph. Under some conditions like (h1)-(h3) on $f$ and appropriate assumptions on $a(x)$, by the mountain pass theorem, they proved that the above-mentioned equation has a positive solution and then by the method of Nehari manifold, they also obtained the existence of  ground state solution.
\par
Inspired by the above papers, in this paper, we investigate the multiplicity of solutions for system (\ref{aa1}) on finite graphs.
It is easy to see that all of these conditions (f1)-(f2) and (h1)-(h3) imply that the nonlinear term $f(x,t)$ needs to have growth both near $0$ and $\infty$ about $t$.
In the present paper, by applying the cut-off technique in \cite{Costa-Wang2005}, we only need to make some assumptions on the nonlinear term $F(x,u,v)$ in a neighborhood of $(u,v)=(0,0)$.
For the super-linear case, by means of the mountain pass theorem, we show that system (\ref{aa1}) has at least a nontrivial weak solution for all $\lambda$ sufficiently large.
Moreover, motivated by \cite{Kang-Liu2020}, we manage to find out the specific lower bound of parameter and the tendency of solution with the change of parameter.
Interestingly, it is necessary to take into account the specific structure of the graph when calculating the specific lower bound of parameter, which is different from the Euclidean setting (see Remark \ref{remark5.1} below).
For the sub-linear case, by using a revised Clark's theorem, we obtain that system (\ref{aa1}) has a sequence of solutions tending to 0 for every $\lambda>0$. In details, we obtain the following results.

\vskip2mm
\par
{\bf(I) The super-linear case:}
\vskip2mm
\par
\noindent
\begin{theorem}\label{theorem1.1}
Let $G=(V,E)$ be a finite graph, $m_{i}\geq 1$ and $p_{i}> 1$, $i=1,2$, $\lambda$ is a parameter with $\lambda>0$, $h_{i}(x)\in C(V,\R ^{+})$, $i=1,2$, $F:V\times\R\times\R \rightarrow \R$ satisfies the following conditions:
\begin{itemize}
\item[$(H_0)$ ] there exists a constant $\delta>0$ such that $F(x,t,s)$ is continuously differentiable in $(t,s)\in \R\times\R$ with $|(t,s)|\leq\delta$ for all $x\in V$;
\item[$(H_1)$ ] there exist $q_{i}>p_{i}$ and $M_{i}>0$, $i=1,2$ such that
 $$  F(x,t,s)\geq M_{1}|t|^{q_{1}}+M_{2}|s|^{q_{2}}$$
 for all $x\in V$ and $|(t,s)|\leq\delta$;
\item[$(H_2)$ ] there exist $k_{i}\in(p_{i},q_{i}) (i=1,2)$
satisfy $\min\{k_{1},k_{2}\}>\max\{p_{1},p_{2}\}$
and $M_{j}>0(j=3,4)$ such that
$$  |F_{t}(x,t,s)|\leq M_{3}|t|^{k_{1}-1}+M_{4}|s|^{\frac{k_{2}(k_{1}-1)}{k_{1}}},$$
$$  |F_{s}(x,t,s)|\leq M_{3}|t|^{\frac{k_{1}(k_{2}-1)}{k_{2}}}+M_{4}|s|^{k_{2}-1}$$
for all $x\in V$ and $|(t,s)|\leq\delta$, where $F_{t}(x,t,s)=\frac{\partial F(x,t,s)}{\partial t}$ and  $F_{s}(x,t,s)=\frac{\partial F(x,t,s)}{\partial s}$;
\item[$(H_3)$ ]there exist constants $\beta_{i}>p_{i}$, $i=1,2$ such that
$$0<F(x,t,s)\leq \frac{1}{\beta_{1}}tF_{t}(x,t,s)+ \frac{1}{\beta_{2}}sF_{s}(x,t,s)$$
for all $x\in V$ and $|(t,s)|\leq\delta$ with $(t,s)\neq (0,0)$.
\end{itemize}
Then there exists at least a nontrivial weak solution $(u_{\lambda},v_{\lambda})$ to system (\ref{aa1})
for all
$\lambda>\lambda_{*}:=\max\{\Lambda_1,\Lambda_2,\Lambda_3,\Lambda_4,\Lambda_5\}$ and
$\lim_{\lambda\to \infty} \|(u_{\lambda},v_{\lambda})\|=0=\lim_{\lambda\to \infty} \|(u_{\lambda},v_{\lambda})\|_\infty$,
where
$\theta_{i}=\min\{\beta_{i},k_{i}\}$, $i=1,2$,
$(u_{0},v_{0})\in W^{m_{1},p_{1}}(V)\times W^{m_{2},p_{2}}(V)$
satisfies
$0<\|u_{0}\|_{\infty}\leq \frac{\delta}{2}$,
$0<\|v_{0}\|_{\infty}\leq \frac{\delta}{2}$
for any $x\in V$,
$M_{5}=\frac{1}{k_{1}}(M_{3}+M_{4})$,
$M_{6}=(\frac{k_{1}-1}{k_{1}}+\frac{1}{k_{2}})M_{4}$,
$h_{i,\infty}=\min\limits_{x\in V} h_{i}(x)$,
$h_{i}^{\infty}=\max\limits_{x\in V} h_{i}(x)$, $i=1,2$,
$\mu_\infty=\min\limits_{x\in V} \mu(x)$,
 \begin{eqnarray}
 & &
       \|(u_{\lambda},v_{\lambda})\|
 =     \|u_{\lambda}\|_{m_{1},p_{1}}+\|v_{\lambda}\|_{m_{2},p_{2}},
 \quad
       \|(u_{\lambda},v_{\lambda})\|_\infty
 =
       \max_{x\in V}\ |(u_{\lambda}(x),v_{\lambda}(x))|,
 \nonumber\\
 & &
 \label{A1}
        \Lambda_1
 =
       \frac{2^{1-\max\{p_{1},p_{2}\}}}
      {
      \max\{p_{1},p_{2}\}
      }
     \left(
     \frac{M_{5}}{\mu_{\infty}^{\frac{k_{1}-p_{1}}{p_{1}}} h_{1,\infty}^{\frac{k_{1}}{p_{1}}}}
      +
     \frac{M_{6}}{\mu_{\infty}^{\frac{k_{2}-p_{2}}{p_{2}}} h_{2,\infty}^{\frac{k_{2}}{p_{2}}}}
      \right)^{-1},
\\
 & &
 \label{A2}
        \Lambda_2
=
   \frac{2^{1-\max\{p_{1},p_{2}\}}}{\max\{p_{1},p_{2}\}}
   \left(
  \frac{M_{5}}{\mu_{\infty}^{\frac{k_{1}-p_{1}}{p_{1}}} h_{1,\infty}^{\frac{k_{1}}{p_{1}}}}
  +
  \frac{M_{6}}{\mu_{\infty}^{\frac{k_{2}-p_{2}}{p_{2}}} h_{2,\infty}^{\frac{k_{2}}{p_{2}}}}
  \right)^{-1}
  \left(
   6h_{1,\infty}^{\frac{1}{p_{1}}}\|u_{0}\|_{L^{p_{1}}(V)}
  \right)^{\max\{p_{1},p_{2}\}-\min\{k_{1},k_{2}\}},
\\
 & &
 \label{A3}
        \Lambda_3
=
    \frac{2^{1-\max\{p_{1},p_{2}\}}}{\max\{p_{1},p_{2}\}}
  \left(
  \frac{M_{5}}{\mu_{\infty}^{\frac{k_{1}-p_{1}}{p_{1}}} h_{1,\infty}^{\frac{k_{1}}{p_{1}}}}
  +
  \frac{M_{6}}{\mu_{\infty}^{\frac{k_{2}-p_{2}}{p_{2}}} h_{2,\infty}^{\frac{k_{2}}{p_{2}}}}
  \right)^{-1}
  \left(
  3h_{2,\infty}^{\frac{1}{p_{2}}}\|v_{0}\|_{L^{p_{2}}(V)}
  \right)^{\max\{p_{1},p_{2}\}-\min\{k_{1},k_{2}\}},
 \end{eqnarray}
\begin{eqnarray}
\label{A4}
     \Lambda_4
=
     \dfrac{
      \frac{1}{p_{1}}\max\{1,h_{1}^{\infty}\}
      \left(
      \|\nabla^{m_{1}} u_{0}\|_{L^{p_{1}}(V)}^{p_{1}}
      +\|u_{0}\|_{L^{p_{1}}(V)}^{p_{1}}
      \right)
      +
     \frac{1}{p_{2}}\max\{1,h_{2}^{\infty}\}
     \left(
      \|\nabla^{m_{2}} v_{0}\|_{L^{p_{2}}(V)}^{p_{2}}
      +\|v_{0}\|_{L^{p_{2}}(V)}^{p_{2}}
     \right)
      }
      {
      M_{1}\|u_{0}\|_{L^{p_{1}}(V)}^{q_{1}}
     \left(\sum\limits_{x\in V}\mu(x)\right)
      ^{1-\frac{q_{1}}{p_{1}}}+M_{2}\|v_{0}\|_{L^{p_{2}}(V)}^{q_{2}}
      \left(\sum\limits_{x\in V} \mu(x)\right)
      ^{1-\frac{q_{2}}{p_{2}}}
       },
\end{eqnarray}
\begin{eqnarray}
\label{A5}
         \Lambda_5
=
       \max\Bigg\{
       \left(
       \frac{2^{2p_{1}+1}p_{1}\theta_{1}\max\{C_{1,*},C_{2,*}\}}
       {
       \delta^{p_{1}} (\theta_{1}-p_{1})\mu_{\infty} h_{1,\infty}
       - 2^{2p_{1}+1}p_{1}\theta_{1}\max\{C_{1,*},C_{2,*}\}
       }
       \right)
       ^{\frac{1}{\max\left\{\frac{p_{1}}{q_{1}-p_{1}},\frac{p_{2}}{q_{2}-p_{2}}\right\}}},
\nonumber\\
       \left(
       \frac{2^{2p_{2}+1}p_{2}\theta_{2}\max\{C_{1,*},C_{2,*}\}}
       {
       \delta^{p_{2}} (\theta_{2}-p_{2})\mu_{\infty} h_{2,\infty}
       - 2^{2p_{2}+1}p_{2}\theta_{2}\max\{C_{1,*},C_{2,*}\}
       }
       \right)
       ^{\frac{1}{\max\left\{\frac{p_{1}}{q_{1}-p_{1}},\frac{p_{2}}{q_{2}-p_{2}}\right\}}}
       \Bigg\},
\\
\label{A6}
        C_{1,*}
=
       \frac{(q_{1}-p_{1}) \sum\limits_{x\in V}\mu(x) }{p_{1}}
       \left(
       \frac{(\max\{1,h_{1}^{\infty}\})^{q_{1}}}
       {q_{1} ^{q_{1}} M_{1} ^{p_{1}}}
       \right)
       ^{\frac{1}{q_{1}-p_{1}}}
      \left(
      \frac{ \left(\|\nabla^{m_{1}} u_{0}\|_{L^{p_{1}}(V)}^{p_{1}}
           + \|u_{0}\|_{L^{p_{1}}(V)}^{p_{1}}\right)^{\frac{1}{p_{1}}}}
      { \|u_{0}\|_{L^{p_{1}}(V)} }
      \right)
      ^{\frac{p_{1}q_{1}}{q_{1}-p_{1}}},
\\
\label{A7}
       C_{2,*}
=
       \frac{(q_{2}-p_{2}) \sum\limits_{x\in V}\mu(x) }{p_{2}}
       \left(
       \frac{(\max\{1,h_{2}^{\infty}\})^{q_{2}}}
       {q_{2} ^{q_{2}} M_{2} ^{p_{2}}}
       \right)
       ^{\frac{1}{q_{2}-p_{2}}}
       \left(
       \frac{
       \left(\|\nabla^{m_{2}} v_{0}\|_{L^{p_{2}}(V)}^{p_{2}}
       + \|v_{0}\|_{L^{p_{2}}(V)}^{p_{2}}
       \right)
       ^{\frac{1}{p_{2}}}
       }
       {
       \|v_{0}\|_{L^{p_{2}}(V)}
       }
       \right)
       ^{\frac{p_{2}q_{2}}{q_{2}-p_{2}}}.
\end{eqnarray}
\end{theorem}
\par
\noindent
\begin{remark}\label{remark1.1}
It follows from $(H_1)$ and $(H_3)$ that
\begin{itemize}
\item[$(H_4)$] $F(x,0,0)=0$ for all $x\in V$.
\end{itemize}
Then by $(H_1)$, $(H_2)$ and $(H_4)$, it is easy to check that the following condition:
\begin{itemize}
\item[$(H_5)$] $$M_{1}|t|^{q_{1}}+M_{2}|s|^{q_{2}} \leq F(x,t,s)\leq M_{5}|t|^{k_{1}}+M_{6}|s|^{k_{2}}$$
for all $x\in V$ and $|(t,s)|\leq\delta$, where $M_{5}=\frac{1}{k_{1}}(M_{3}+M_{4})$, $M_{6}=(\frac{k_{1}-1}{k_{1}}+\frac{1}{k_{2}})M_{4}$.
\end{itemize}
 \end{remark}

\par
{\bf(II) The sub-linear case:}
\par
\noindent
\begin{theorem}\label{theorem1.2}
Let $G=(V,E)$ be a finite graph, $m_{i}\geq 1$, $p_{i}> 1$, $i=1,2$, $\lambda$ is a parameter with $\lambda>0$, $h_{i}(x)\in C(V,\R ^{+})$, $i=1,2$, $F:V\times\R\times\R \rightarrow \R$ satisfies $(H_0)$, $(H_4)$ and the following conditions:
\begin{itemize}
\item[$(R_1)$ ] there exist $q_{i}\in(1,p_{i}) (i=1,2)$ satisfy $\min\{p_{1},p_{2}\}>\max\{q_{1},q_{2}\}$ and $K_{i}>0 (i=1,2)$ such that
 $$  F(x,t,s)\geq K_{1}|t|^{q_{1}}+K_{2}|s|^{q_{2}}$$
 for all $x\in V$ and $|(t,s)|\leq\delta$;
\item[$(R_2)$ ] there exist $k_{i}\in(1,q_{i}) (i=1,2)$ and $K_{j}>0(j=3,4)$ such that
$$  |F_{t}(x,t,s)|\leq K_{3}|t|^{k_{1}-1}+K_{4}|s|^{\frac{k_{2}(k_{1}-1)}{k_{1}}},$$
$$  |F_{s}(x,t,s)|\leq K_{3}|t|^{\frac{k_{1}(k_{2}-1)}{k_{2}}}+K_{4}|s|^{k_{2}-1},$$
for all $x\in V$ and $|(t,s)|\leq\delta$;
\item[$(R_3)$ ] $F(x,t,s)=F(x,-t,-s)$ for all $x\in V$ and $|(t,s)|\leq\delta$.
\end{itemize}
Then for every $\lambda>0$,  system (\ref{aa1}) has a sequence of weak solutions $\{(u_{k}^{\lambda},v_{k}^{\lambda})\}$ with $\|(u_{k}^{\lambda},v_{k}^{\lambda})\|\rightarrow 0$ as $k\rightarrow\infty$.
\end{theorem}
\par
\noindent
\begin{remark}\label{remark1.3}
By $(R_1)$, $(R_2)$ and $(H_4)$, it is easy to check that the following condition:
\begin{itemize}
\item[$(R_4)$]
$$
   K_{1}|t|^{q_{1}}+K_{2}|s|^{q_{2}} \leq F(x,t,s)\leq K_{5}|t|^{k_{1}}+K_{6}|s|^{k_{2}}
$$
for all $x\in V$, $(t,s)\in \R^{2}$ with $|(t,s)|\leq\delta$ , where $K_{5}=\frac{1}{k_{1}}(K_{3}+K_{4})$, $K_{6}=(\frac{k_{1}-1}{k_{1}}+\frac{1}{k_{2}})K_{4}$.
\end{itemize}
\end{remark}
\vskip2mm
\par
We organize this paper as follows. In Section 2, we mainly recall some knowledge for poly-Laplacian on graphs and Sobolev spaces. In Section 3, we complete the proof of Theorem \ref{theorem1.1}.
In Section 4, we complete the proof of Theorem \ref{theorem1.2}.
In Section 5, we apply Theorem \ref{theorem1.1} to an example and compute the value of lower bound $\lambda_{*}$. Moreover, we also present one example to illustrate Theorem \ref{theorem1.2}.
\vskip2mm
 {\section{Preliminaries}}
 \vskip2mm
 \noindent
 \par
 In this section, we recall some notions and important properties about poly-Laplacian on graphs and Sobolev spaces and some useful lemmas. One can see details in \cite{Grigor-Lin2016,Han-Shao2021,Zhang-Zhao2018}.
 \par
 First of all, we recall the basic notions of weighted graphs. Let $G=(V,E)$ be a finite graph, where $V$ denotes the vertex set and $E$ denotes the edge set. We say that a graph $G=(V,E)$ is finite if and only if $V$ and $E$ are both finite sets. Fix an edge weight function $\omega:E\rightarrow (0,+\infty)$ satisfying $\omega_{xy}=\omega_{yx}$, for any edge $\{x,y\}\in E$. Let $deg(x)=\Sigma_{x\sim y}\omega_{xy}$ be the degree of $x\in V$, where we write $x\sim y$ if $ \{x,y\}\in E$. Let $\mu:V\rightarrow (0,+\infty)$ be a finite measure on $V$.
 \par
 Define $C(V)$ as the set of all real functions on $V$. Then $C(V)$ is a finitely dimensional linear space with the usual functions additions and scalar multiplications. For any $\varphi \in C(V)$, we denote
\begin{eqnarray}\label{b1}
\int_{V} \varphi d\mu = \sum_{x\in V} \mu(x) \varphi(x).
\end{eqnarray}
First of all, for $\varphi_{1},\varphi_{2} \in C(V)$, the associated gradient form reads
 \begin{eqnarray}\label{b2}
  \Gamma(\varphi_{1},\varphi_{2})(x) = \frac{1}{2\mu(x)} \sum_{y\sim x} \omega_{xy} (\varphi_{1}(y)-\varphi_{1}(x))(\varphi_{2}(y)-\varphi_{2}(x)).
 \end{eqnarray}
Sometimes we use the notation $\nabla \varphi_{1}\nabla \varphi_{2}=\Gamma(\varphi_{1},\varphi_{2})$. For any $\varphi \in C(V)$, we write $\Gamma(\varphi) = \Gamma(\varphi,\varphi)$. Then the length of its gradient is defined as
\begin{eqnarray}\label{b3}
|\nabla \varphi|
=  \sqrt{\Gamma(\varphi)}
=  \left(\frac{1}{2\mu(x)} \sum_{y\sim x} \omega_{xy} (\varphi(y)-\varphi(x))^{2}\right)^{\frac{1}{2}}.
  \end{eqnarray}
We have the following equality:
\begin{eqnarray*}
\int_{V} (\Delta \varphi)\phi d\mu = -\int_{V}  \Gamma(\varphi,\phi)d\mu, \ \ \forall \phi \in C_{c}(V),
  \end{eqnarray*}
where the general discrete graph Laplacian $ \Delta:C(V)\rightarrow C(V) $ is defined as
 \begin{eqnarray*}
  \Delta \varphi(x) = \frac{1}{\mu(x)} \sum_{y\sim x} \omega_{xy} (\varphi(y)-\varphi(x)).
 \end{eqnarray*}
By Lemma 2.1 in \cite{Zhang-Zhao2018}, we know that $ \Delta \varphi$ is well-defined.
Next, the length of $m$-order gradient of $\varphi \in C(V)$ is defined by
\begin{eqnarray*}
|\nabla^{m}\varphi|= \begin{cases}
               |\Delta^{\frac{m}{2}} \varphi|, & \mbox{  when } m\  \mbox{ is even},\\
               |\nabla \Delta^{\frac{m-1}{2}} \varphi |,&  \mbox{ when } m\  \mbox{ is odd} ,
\end{cases}
\end{eqnarray*}
where $m\geq 1$, $|\Delta^{\frac{m}{2}} \varphi| $ is defined as the usual absolute of $\Delta^{\frac{m}{2}} \varphi $ and $|\nabla \Delta^{\frac{m-1}{2}} \varphi|$ denotes (\ref{b3}) with $\varphi$ replaced by $\nabla \Delta^{\frac{m-1}{2}} \varphi$.
For any $p> 1$, the $p$-Laplacian operator of $\varphi \in C(V)$, denoted by $\Delta_{p}\varphi $, is defined as the following form
\begin{eqnarray*}
\int_{V} (\Delta_{p} \varphi)\phi d\mu = -\int_{V} |\nabla \varphi|^{p-2} \Gamma(\varphi,\phi)d\mu, \ \ \forall \phi \in C_{c}(V),
\end{eqnarray*}
where $ \Gamma(\varphi,\phi)$ is defined as in (\ref{b2}), $C_{c}(V)$ denotes the set of all real functions with compact support, and the integration is defined by (\ref{b1}). $ \Delta_{p}:C(V)\rightarrow C(V) $ can also be written as
\begin{equation}\label{c1}
  \Delta_{p}\varphi(x) = \frac{1}{2\mu(x)} \sum_{y\sim x}  \left(|\nabla \varphi|^{p-2} (y)+ |\nabla \varphi|^{p-2} (x) \right) \omega_{xy} (\varphi(y)-\varphi(x)).
  \end{equation}
By Lemma 2.1 in \cite{Han-Shao2021}, we know that the definition of $p$-Laplacian operator in (\ref{c1}) is reasonable. Now, the generalized graph poly-Laplacian $\pounds_{m,p}:C(V)\rightarrow C(V) $ (see \cite{Grigor-Lin2016}) is introduced as  following
\begin{equation*}
\int_{V} (\pounds_{m,p}\varphi)\phi d\mu =
 \begin{cases}
 \int_{V}|\nabla^{m} \varphi|^{p-2} \Gamma(\Delta^{\frac{m-1}{2}} \varphi, \Delta^{\frac{m-1}{2}} \phi)d\mu, \;\;\; \mbox{ when } m\  \mbox{ is odd},\\
  \int_{V}|\nabla^{m}\varphi|^{p-2} \Delta^{\frac{m}{2}}\varphi \Delta^{\frac{m}{2}} \phi d\mu ,\;\;\; \mbox{ when } m\  \mbox{ is even},
   \end{cases}
 \end{equation*}
where $m\geq1$, $p>1$, $\varphi \in C(V)$, $\phi \in C_{c}(V)$. It is easy to see that $\pounds_{m,2}=(-\Delta)^{m}$ and $\pounds_{1,p}=-\Delta_{p}$.
 \par
 For any $1\leq q <\infty$, assume that the completion space of $C_{c}(V)$ is $L^{q}(V)$ under the norm
 $$\|\varphi \|_{L^{q}(V)}=\left(\int_{V} |\varphi|^{q} d\mu \right)^{\frac{1}{q}}.  $$
 Moreover, the completion space of $C_{c}(V)$ is $W^{m,p}(V)$ under the norm
 \begin{eqnarray}\label{b7}
  \|\varphi\|_{m,p}=\left(\int_V ( |\nabla^{m} \varphi|^p+h(x)|\varphi|^p)d\mu\right)^{1/p},
  \end{eqnarray}
 where $m\geq1$, $p>1$, $h(x) >0$ for all $x\in V$. Observing that $W^{m,p}(V)$ is a finite dimensional linear space and $W^{m,p}(V)$ with norm (\ref{b7}) is a Banach space. For more details about poly-Laplacian operator and Sobolev spaces on graphs, we refer readers to \cite{Grigor-Lin2016} and the reference therein.

\par
The following inequalities and embedding proposition involving  Sobolev spaces will be used frequently in our proofs.
\par
\noindent
\begin{remark}\label{remark2.1}
It is easy to obtain that
  $$
  h_{\infty}\|\varphi\|_{L^{p}(V)}^{p} \leq \|\varphi\|_{m,p}^{p} \leq  \max\{1,h^{\infty}\}\left(\|\nabla^m \varphi\|_{L^{p}(V)}^{p}+\|\varphi\|_{L^{p}(V)}^{p}\right),
  $$
  where $h_\infty=\min\limits_{x\in V} h(x)$ and $h^\infty=\max\limits_{x\in V} h(x)$.
\end{remark}
\par
\noindent
\begin{lemma}(\cite{Zhang-Zhang-Xie-Yu-2022})\label{lemma2.1}
Let $p>1$. For all $\varphi\in W^{m,p}(V)$, it holds
  \begin{equation*}
   \|\varphi\|_{\infty}\leq \frac{1}{\mu_\infty^{\frac{1}{p}} h_\infty^{\frac{1}{p}}}\|\varphi\|_{m,p},
  \end{equation*}
  where $\mu_\infty=\min\limits_{x\in V} \mu(x)$.
\end{lemma}
\par
\noindent
\begin{lemma}(\cite{Zhang-Zhang-Xie-Yu-2022})\label{lemma2.2}
Let $G=(V,E)$ be a finite graph.
Let $m$ be any positive integer and $p>1$.
Then the embedding $W^{m,p}(V)\hookrightarrow L^{q}(V)$ is continuous for all $1\leq q<+\infty$ and
\begin{eqnarray*}
\|\varphi\|_{L^q(V)}\leq \frac{\left(\sum_{x\in V}\mu(x)\right)^{1/q}}{\mu_\infty^{\frac{1}{p}} h_\infty^{\frac{1}{p}}}\|\varphi\|_{m,p}
\end{eqnarray*}
for all $\varphi\in W^{m,p}(V)$. Moreover, $W^{m,p}(V)$ is pre-compact.
\end{lemma}
\par
Next, we give an inequality which will be used in the proof of Theorem \ref{theorem1.1}.
\par
\noindent
\begin{lemma}\label{lemma2.3}
Let $\lambda>0$, $a>0$ and $b>0$. Then
$$
\lambda^{-a}+\lambda^{-b}\leq 2\left( 1+ \lambda^{-\max\{a ,b \}}   \right).
$$
\end{lemma}
\vskip2mm
 \noindent
 {\bf Proof.} If $0<\lambda<1$, we have $\lambda^{-a}+\lambda^{-b}\leq2\lambda^{-\max\{a ,b \}}$. If $\lambda\geq 1$, it is obvious that $\lambda^{-a}+\lambda^{-b}\leq2$. Then for any $\lambda>0$, the conclusion holds.
\qed

\vskip2mm
\par
Now, we recall a special version of the mountain pass lemma. Let $X$ be a real Banach space, $I \in C^{1}(X,\R)$ and $c\in \R$. A sequence $\{\varphi_{n}\}\subset X$ is said to be a $(C)_{c}$-sequence of $I$ in $X$ if
 \begin{equation*}
   I(\varphi_{n})\to c,\,\ \|I'(\varphi_{n})\|_{X^{*}}(1+\|\varphi_{n}\|)\to 0 .
   \end{equation*}
$I$ is said to satisfy the $(C)_{c}$-condition if any $(C)_{c}$-sequence has a convergent subsequence. To prove Theorem \ref{theorem1.1}, we need the following mountain pass theorem.
\par
\noindent
\begin{lemma}(\cite{Ekeland1990})\label{lemma2.4}
Let $X$ be a real Banach space, $I \in C^{1}(X,\R)$, $w\in X$ and $r>0$ be such that $\|w\|>r$ and
 $$
 \inf_{\|\varphi\|=r} I(\varphi)>I(0)=0 \geq I(w).
 $$
 Then there exists a $(C)_{c}$-sequence with
 $$
 c:=\inf_{\gamma\in\Gamma}\max_{t\in[0,1]}I(\gamma(t)),
 $$
 $$
 \Gamma:=\{\gamma\in C([0,1],X):\gamma(0)=0, \; \gamma(1)= w\}.
 $$
\end{lemma}
\par
Next, we recall an extension of Clark's theorem. Let $X$ be a real Banach space and $I \in C^{1}(X,\R)$. For any sequence $\{\varphi_{n}\}\subset X$, if $\{I(\varphi_{n})\}$ is bounded and $I'(\varphi_{n})\to 0$ as $n\to \infty$, then $\{\varphi_{n}\}$ is said to be a Palais-Smale sequence of $I$ in $X$. If any Palais-Smale sequence $\{\varphi_{n}\}$ admits a convergent subsequence, then we call that $I$  satisfies Palais-Smale ((PS) for short) condition . To prove Theorem \ref{theorem1.2}, we need the following Clark's theorem.
\par
\noindent
\begin{lemma}(\cite{Liu-Wang2015})\label{lemma2.5}
Let $X$ be a real Banach space and $I \in C^{1}(X,\R)$. Suppose that $I$ satisfies the (PS)-condition, which is even, bounded from below, and $I(0)=0$. If for any $k\in \mathbb{N}$, there exists a $k$-dimensional subspace $Z^{k}$ of $X$ and $\rho_{k}>0$ such that $\sup_{Z^{k}\cap S_{\rho_{k}}} I<0$, where $S_{\rho_{k}}=\{\varphi\in X|\|\varphi\|=\rho_{k}  \}$, then at least one of the following conclusions holds:
\par
 (i) There exist a sequence of critical points $\{\varphi_{k}\}$ satisfying $I(\varphi_{k}) < 0$ for all $k$ and $\|\varphi_{k}\|\rightarrow 0$ as $k \rightarrow \infty$.
 \par
 (ii) There exist a constant $\zeta > 0$ such that for any $0<a<\zeta$ there exists a critical point $\varphi$ such that $\|\varphi\|= a$ and $I(\varphi)=0$.
\end{lemma}
\par
\noindent
\begin{remark}(\cite{Liu-Wang2015})\label{remark2.6}
Lemma \ref{lemma2.5} implies that there exist a sequence of critical points $\varphi_{k}\neq0$ such that $I(\varphi_{k})\leq0$, $I(\varphi_{k})\rightarrow 0$ and $\|\varphi_{k}\|\rightarrow 0$ as $k \rightarrow \infty$.
\end{remark}
\vskip2mm
 {\section{Proofs for the super-linear case}}
  \setcounter{equation}{0}
\par
We define the space $W := W^{m_{1},p_{1}}(V)\times W^{m_{2},p_{2}}(V)$ with the norm $\|(u,v)\|=\|u\|_{m_{1},p_{1}}+ \|v\|_{m_{2},p_{2}}$, where $m_{i}\geq1$, $p_{i}>1$, $i=1,2$. Then $(W, \|\cdot\|) $ is a Banach space. On $W$, we define the variational functional corresponding to system (\ref{aa1}) by
 \begin{eqnarray*}
        I_{\lambda}(u,v)
    =   \frac{1}{p_{1}}\int_{V}(|\nabla^{m_{1}} u|^{p_{1}}
     +h_{1}(x)|u|^{p_{1}})d\mu + \frac{1}{p_{2}}\int_{V}(|\nabla^{m_{2}}v|^{p_{2}}
     +h_{2}(x)|v|^{p_{2}})d\mu
      -\lambda \int_V F(x,u,v)d\mu.
  \end{eqnarray*}
Since conditions $(H_1)$ and $(H_2)$ imply the behavior of $F$ just near the origin, the functional $\int_V F(x,u,v)d\mu$ is not well defined in the Sobolev space $W$.
To deal with this problem, we need to extend $F$ to a proper function $\bar{F}$ by the cut-off technique developed by \cite{Costa-Wang2005}.
In order to adapt system (\ref{aa1}), we make an extension to $\R ^{2}$ for the cut-off function in \cite{Costa-Wang2005}.
\par
Define $\tau(t,s)\in C^1(\R ^{2},[0,1])$ as a cut-off function fulfilling $t\tau'_{t}(t,s)\leq 0$ and $s\tau'_{s}(t,s)\leq 0$ for all $(t,s)\in\R ^{2}$ with $\frac{\delta}{2}<|(t,s)|\leq\delta$ and
 \begin{eqnarray*}
 \tau(t,s)= \begin{cases}
           1, \;\;\;  \text { if }\;\;|(t,s)| \leq \delta/2, \\
           0, \;\;\;  \text { if }\;\;|(t,s)|    >   \delta.
          \end{cases}
   \end{eqnarray*}
Define $\bar{F}:V\times\R\times \R\to\R$ by
   $$
   \bar{F}(x,t,s)=\tau(t,s)F(x,t,s)+(1-\tau(t,s))\left(M_{5}|t|^{k_{1}}+ M_{6}|s|^{k_{2}}\right).
   $$
By the definition of $\tau(t,s)$, we could obtain the following lemma:
\par
\noindent
\begin{lemma}\label{lemma3.1}
Assume that $(H_0)$, $(H_2)$, $(H_3)$ and $(H_5)$ hold. Then
\begin{itemize}
\item[$(H_0)'$ ]$\bar{F}(x,t,s)$ is continuously differentiable in $(t,s)\in\R\times \R$ for all $x \in V$;
\item[$(H_2)'$ ] there exist $M_{7}>0$ and $M_{8}>0$ such that
$$  |\bar{F}_{t}(x,t,s)|\leq M_{7}|t|^{k_{1}-1}+M_{8}|s|^{\frac{k_{2}(k_{1}-1)}{k_{1}}},$$
$$  |\bar{F}_{s}(x,t,s)|\leq M_{7}|t|^{\frac{k_{1}(k_{2}-1)}{k_{2}}}+M_{8}|s|^{k_{2}-1}$$
for all $x\in V$ and $(t,s)\in \R ^{2}$;
\item[$(H_3)'$ ]$$0<\bar{F}(x,t,s)\leq \frac{1}{\theta_{1}}t\bar{F}_{t}(x,t,s)+ \frac{1}{\theta_{2}}s\bar{F}_{s}(x,t,s)$$
for all $x\in V$ and  $(t,s)\in \R ^{2}\backslash \{(0,0)\}$, where $\theta_{i}=\min\{\beta_{i},k_{i}\}$, $i=1,2$;
\item[$(H_5)'$ ]$$
   M_{1}|t|^{q_{1}}+M_{2}|s|^{q_{2}} \leq \bar{F}(x,t,s)\leq M_{5}|t|^{k_{1}}+M_{6}|s|^{k_{2}}
 $$
 for all $x\in V$ and $|(t,s)|\leq\delta$,
 $$
  0 \leq \bar{F}(x,t,s)\leq M_{5}|t|^{k_{1}}+M_{6}|s|^{k_{2}}
 $$
 for all $x\in V$ and $(t,s)\in \R ^{2}$.
\end{itemize}
\end{lemma}
\par
Consider the modified system of (\ref{aa1}) given by
\begin{equation}\label{mod1}
 \begin{cases}
 \pounds_{m_{1},p_{1}}u+h_{1}(x)|u|^{p_{1}-2}u=\lambda \bar{F}_{u}(x,u,v),  \;\;\; x \in V,\\
 \pounds_{m_{2},p_{2}}v+h_{2}(x)|v|^{p_{2}-2}v=\lambda \bar{F}_{v}(x,u,v),  \;\;\; x \in V.
   \end{cases}
 \end{equation}
 Define the corresponding functional $\bar{I}_{\lambda}: W\rightarrow \R$ by
 \begin{eqnarray*}
        \bar{I}_{\lambda}(u,v)
  &  : =  & \frac{1}{p_{1}}\int_{V}(|\nabla^{m_{1}} u|^{p_{1}}
     +h_{1}(x)|u|^{p_{1}})d\mu + \frac{1}{p_{2}}\int_{V}(|\nabla^{m_{2}}v|^{p_{2}}
     +h_{2}(x)|v|^{p_{2}})d\mu \\
    & &  -\lambda \int_V \bar{F}(x,u,v)d\mu
  \end{eqnarray*}
for all $ (u, v)\in W $. By $(H_0)'$, $(H_2)'$, $(H_5)'$ and the continuous embeddings
 \begin{equation*}
   W^{m_{1},p_{1}}(V)\hookrightarrow L^{q_{1}}(V),\,\,\,\ W^{m_{2},p_{2}}(V)\hookrightarrow L^{q_{2}}(V),
   \end{equation*}
standard arguments show that $\bar{I}_{\lambda}$ is well defined and of class $C^{1}$ on $W$, and
 \begin{eqnarray*}
          \langle \bar{I}_{\lambda}'(u,v),(\tilde{u},\tilde{v})\rangle
  &  =  & \int_V\left(\pounds_{m_{1},p_{1}} u \tilde{u}+h_{1}(x)|u|^{p_{1}-2}u\tilde{u}\right)d\mu\\
  &    & + \int_V\left(\pounds_{m_{2},p_{2}} v \tilde{v}+h_{2}(x)|v|^{p_{2}-2}v\tilde{v}\right)d\mu\\
   &    &-\lambda \int_V \bar{F}_{u}(x,u,v)\tilde{u}d\mu- \lambda \int_V \bar{F}_{v}(x,u,v)\tilde{v}d\mu
   \end{eqnarray*}
   for all $ (u,v),(\tilde{u},\tilde{v})\in W$ (for example, see \cite{Yang-Zhang-2024}). Hence
   \begin{eqnarray*}
          \langle \bar{I}_{\lambda}'(u,v),(u,v)\rangle
    =   \|u\|_{m_{1},p_{1}}^{p_{1}}+\|v\|_{m_{2},p_{2}}^{p_{2}}-\lambda \int_V \bar{F}_{u}(x,u,v)ud\mu- \lambda \int_V \bar{F}_{v}(x,u,v)vd\mu
   \end{eqnarray*}
 for all $ (u,v)\in W$. Note that the critical points of $\bar{I}_{\lambda}$ with $ L^{\infty}$-norm less than or equal to $ \delta/2$ are also critical points of $I_{\lambda}$. Then the weak solutions of system (\ref{mod1}) with $ L^{\infty}$-norm less than or equal to $ \delta/2$ are also weak solutions of problem (\ref{aa1}).
 \par
\noindent
\begin{lemma}\label{lemma3.2}
Assume that $(H_5)'$ holds. Then for each $\lambda>0$, there exists two positive constants $\nu_{\lambda}$ and $\eta_{\lambda}$ such that $\bar{I}_{\lambda}(u,v)\geq \eta_{\lambda}$ for all $\|(u,v)\|=\nu_{\lambda}$.
\end{lemma}
 \vskip0mm
   \noindent
 {\bf Proof.}
 For any given $\lambda>0$, $(u,v)\in W$ with $\|(u,v)\|\leq 1$, by $(H_5)'$, Lemma \ref{lemma2.1}, Remark \ref{remark2.1} and $\min\{k_{1},k_{2}\}>\max\{p_{1},p_{2}\}$, we have
\begin{eqnarray}\label{3.2.0}
         \bar{I}_{\lambda}(u,v)
 & = &
        \frac{1}{p_{1}}\|u\|_{m_{1},p_{1}}^{p_{1}}
        +\frac{1}{p_{2}}\|v\|_{m_{2},p_{2}}^{p_{2}}
        -\lambda \int_V \bar{F}(x,u,v)d\mu
\nonumber \\
 &\geq &
        \frac{1}{\max\{p_{1},p_{2}\}}
        (\|u\|_{m_{1},p_{1}}^{\max\{p_{1},p_{2}\}}+\|v\|_{m_{2},p_{2}}^{\max\{p_{1},p_{2}\}})
        - \lambda M_{5}\|u\|_{\infty} ^{k_{1}-p_{1}}\|u\|_{L^{p_{1}}(V)}^{p_{1}}
        -\lambda M_{6}\|v\|_{\infty} ^{k_{2}-p_{2}}\|v\|_{L^{p_{2}}(V)}^{p_{2}}
\nonumber \\
&\geq &
         \frac{ 2^{1-\max\{p_{1},p_{2}\}}}{\max\{p_{1},p_{2}\}}
          \|(u,v)\|^{\max\{p_{1},p_{2}\}}
          -\frac{\lambda M_{5}}{(\mu_{\infty}h_{1,\infty})^{\frac{k_{1}-p_{1}}{p_{1}}}} \|u\|_{m_{1},p_{1}} ^{k_{1}-p_{1}}\|u\|_{L^{p_{1}}(V)}^{p_{1}}
\nonumber \\
&  &
         -\frac{\lambda M_{6}}{(\mu_{\infty},h_{2,\infty})^{\frac{k_{2}-p_{2}}{p_{2}}}} \|v\|_{m_{2},p_{2}} ^{k_{2}-p_{2}}\|v\|_{L^{p_{2}}(V)}^{p_{2}}
\nonumber \\
&\geq &
          \frac{ 2^{1-\max\{p_{1},p_{2}\}}}{\max\{p_{1},p_{2}\}}
          \|(u,v)\|^{\max\{p_{1},p_{2}\}}
          -\frac{\lambda M_{5}}{\mu_{\infty}^{\frac{k_{1}-p_{1}}{p_{1}}} h_{1,\infty}^{\frac{k_{1}}{p_{1}}}} \|u\|_{m_{1},p_{1}} ^{k_{1}}
          -\frac{\lambda M_{6}}{\mu_{\infty}^{\frac{k_{2}-p_{2}}{p_{2}}} h_{2,\infty}^{\frac{k_{2}}{p_{2}}}} \|v\|_{m_{2},p_{2}} ^{k_{2}}
\nonumber \\
&\geq &
          \frac{ 2^{1-\max\{p_{1},p_{2}\}}}{\max\{p_{1},p_{2}\}}
          \|(u,v)\|^{\max\{p_{1},p_{2}\}}
          -\lambda
          \left(
          \frac{M_{5}}{\mu_{\infty}^{\frac{k_{1}-p_{1}}{p_{1}}} h_{1,\infty}^{\frac{k_{1}}{p_{1}}}}
          +
          \frac{M_{6}}{\mu_{\infty}^{\frac{k_{2}-p_{2}}{p_{2}}} h_{2,\infty}^{\frac{k_{2}}{p_{2}}}}
          \right)
           \|(u,v)\|^{\min\{k_{1},k_{2}\}}.
 \end{eqnarray}
Let
\begin{equation*}
  \nu_{0,\lambda}
=
  \left(
  \frac{2^{1-\max\{p_{1},p_{2}\}}}
  {\lambda\max\{p_{1},p_{2}\}
  \left(
  \frac{M_{5}}{\mu_{\infty}^{\frac{k_{1}-p_{1}}{p_{1}}} h_{1,\infty}^{\frac{k_{1}}{p_{1}}}}
  +
  \frac{M_{6}}{\mu_{\infty}^{\frac{k_{2}-p_{2}}{p_{2}}} h_{2,\infty}^{\frac{k_{2}}{p_{2}}}}
  \right)}
  \right)
  ^{\frac{1}{\min\{k_{1},k_{2}\}-\max\{p_{1},p_{2}\}}}.
\end{equation*}
Then by (\ref{3.2.0}), we can choose $\nu_{\lambda}$ satisfies $0<\nu_{\lambda}<\min\{1, \nu_{0,\lambda}\}$ such that
\begin{eqnarray*}
        \bar{I}_{\lambda}(u,v)
&\geq&
       \frac{ 2^{1-\max\{p_{1},p_{2}\}}}{\max\{p_{1},p_{2}\}}
          \nu_{\lambda}^{\max\{p_{1},p_{2}\}}
          -\lambda
          \left(
          \frac{M_{5}}{\mu_{\infty}^{\frac{k_{1}-p_{1}}{p_{1}}} h_{1,\infty}^{\frac{k_{1}}{p_{1}}}}
          +
          \frac{M_{6}}{\mu_{\infty}^{\frac{k_{2}-p_{2}}{p_{2}}} h_{2,\infty}^{\frac{k_{2}}{p_{2}}}}
          \right)
           \nu_{\lambda}^{\min\{k_{1},k_{2}\}}
\nonumber\\
&:=&
         \eta_{\lambda}
\nonumber\\
&>&
         \frac{ 2^{1-\max\{p_{1},p_{2}\}}}{\max\{p_{1},p_{2}\}}
           \nu_{0,\lambda}^{\max\{p_{1},p_{2}\}}
          -\lambda
          \left(
          \frac{M_{5}}{\mu_{\infty}^{\frac{k_{1}-p_{1}}{p_{1}}} h_{1,\infty}^{\frac{k_{1}}{p_{1}}}}
          +
          \frac{M_{6}}{\mu_{\infty}^{\frac{k_{2}-p_{2}}{p_{2}}} h_{2,\infty}^{\frac{k_{2}}{p_{2}}}}
          \right)
            \nu_{0,\lambda}^{\min\{k_{1},k_{2}\}}=0
 \end{eqnarray*}
for all $\|(u,v)\|=\nu_{\lambda}$.
\qed
\par
\noindent
\begin{lemma}\label{lemma3.3}
Assume that $(H_5)'$ holds.
Then for each $\lambda>\max\{\Lambda_1,\Lambda_2,\Lambda_3, \Lambda_4\}$,
where $\Lambda_1$, $\Lambda_2$, $\Lambda_3$ and $\Lambda_4$ are defined in (\ref{A1}), (\ref{A2}), (\ref{A3}) and (\ref{A4}), respectively,
there exists $(u_{0},v_{0})\in W^{m_{1},p_{1}}(V)\backslash\{0\}\times W^{m_{2},p_{2}}(V)\backslash\{0\}$ with $u_{0}>0$, $v_{0}>0$ and $\|(u_{0},v_{0})\|_{\infty}\leq \delta$
such that $\|(u_{0},v_{0})\|>\nu_{\lambda}$ and $\bar{I}_{\lambda}(u_{0},v_{0})<0$.
\end{lemma}
\vskip2mm
   \noindent
 {\bf Proof.}
 Let $(u_{0},v_{0})\in W^{m_{1},p_{1}}(V)\backslash\{0\}\times W^{m_{2},p_{2}}(V)\backslash\{0\}$ with $u_{0}>0$, $v_{0}>0$ for any $x\in V$ and satisfies $\|u_{0}\|_{\infty}\leq \frac{\delta}{2}$ and $\|v_{0}\|_{\infty}\leq \frac{\delta}{2}$.
It is easy to see that
$\nu_{0,\lambda}<1$ for any $\lambda>\Lambda_1$.
We could take  $\nu_{\lambda}=\frac{1}{2}\nu_{0,\lambda}$.
By Remark \ref{remark2.1}, we have $\|u_{0}\|_{m_{1},p_{1}}>\frac{1}{3}\nu_{\lambda}$ for all $\lambda>\Lambda_{2}$
and $\|v_{0}\|_{m_{2},p_{2}}>\frac{2}{3}\nu_{\lambda}$ for all $\lambda>\Lambda_{3}$, respectively.
Then $\|(u_{0},v_{0})\|=\|u_{0}\|_{m_{1},p_{1}}+\|v_{0}\|_{m_{2},p_{2}}
 >\frac{1}{3}\nu_{\lambda}+\frac{2}{3}\nu_{\lambda}=\nu_{\lambda}$  for all $\lambda>\max\{\Lambda_{2},\Lambda_{3}\}$.
It follows from Remark \ref{remark2.1}, $(H_5)'$ and H\"older's inequality that
 \begin{eqnarray*}
         \bar{I}_{\lambda}(u_{0},v_{0})
 & = &
         \frac{1}{p_{1}}\|u_{0}\|_{m_{1},p_{1}}^{p_{1}}
         +\frac{1}{p_{2}}\|v_{0}\|_{m_{2},p_{2}}^{p_{2}}-\lambda \int_V \bar{F}(x,u_{0},v_{0})d\mu
 \\
 &\leq &
       \frac{1}{p_{1}}\max\{1,h_{1}^{\infty}\}(\|\nabla^{m_{1}} u_{0}\|_{L^{p_{1}}(V)}^{p_{1}}+\|u_{0}\|_{L^{p_{1}}(V)}^{p_{1}})
 \\
 & &     +\frac{1}{p_{2}}\max\{1,h_{2}^{\infty}\}(\|\nabla^{m_{2}} v_{0}\|_{L^{p_{2}}(V)}^{p_{2}}+
        \|v_{0}\|_{L^{p_{2}}(V)}^{p_{2}})
 \\
 & &
       -\lambda\left( M_{1}\left(\sum_{x\in V} \mu(x)\right) ^{1-\frac{q_{1}}{p_{1}}}\|u_{0}\|_{L^{p_{1}}(V)}^{q_{1}}+M_{2}\left(\sum_{x\in V} \mu(x)\right) ^{1-\frac{q_{2}}{p_{2}}}\|v_{0}\|_{L^{p_{2}}(V)}^{q_{2}}\right)
 \\
 & <  &  0
 \end{eqnarray*}
for all $\lambda>\Lambda_4$.
\qed
\vskip2mm
 \par
It follows from Lemma \ref{lemma3.2}, Lemma \ref{lemma3.3} and the fact $\bar{I}_{\lambda}(0,0)=0$  that $\bar{I}_{\lambda}$ has a mountain pass geometry, that is, setting
 $$
 \Gamma:=\{\gamma\in C([0,1],W):\gamma(0)=(0,0), \; \gamma(1)=(u_{0},v_{0})\}.
 $$
 Obviously,  $\Gamma\neq \emptyset$.
 By Lemma \ref{lemma2.4}, for the mountain pass level
 $$
 c_{\lambda}:=\inf_{\gamma\in\Gamma}\max_{t\in[0,1]}\bar{I}_{\lambda}(\gamma(t)),
 $$
 there exists a $(C)_{c_{\lambda}}$-sequence $\{(u_{n},v_{n})\}:=\{(u_{n,\lambda},v_{n,\lambda})\}$ of $\bar{I}_{\lambda}$ in $W$, that is,
 \begin{equation}\label{3.2.2}
   \bar{I}_{\lambda}(u_{n},v_{n})\to c_{\lambda}\ \ \mbox{and }\ \ (1+\|(u_{n},v_{n})\|)\|\bar{I}_{\lambda}'(u_{n},v_{n})\|_{W^{*}}\to 0,\ \ \mbox{as } n\to \infty.
   \end{equation}
Moreover, Lemma \ref{lemma3.2} implies that $c_{\lambda}>0$.
\par
\noindent
\begin{lemma}\label{lemma3.4}
Assume that $(H_3)'$ holds. Then the $(C)_{c_{\lambda}}$-sequence $\{(u_{n},v_{n})\}$ has a convergent subsequence for any given $\lambda>0$ .
\end{lemma}
\vskip0mm
   \noindent
   {\bf Proof.} By (\ref{3.2.2}), we have
\begin{eqnarray*}
 \|(u_{n},v_{n})\|\|\bar{I}_{\lambda}'(u_{n},v_{n})\|_{W^{*}}\to 0 \ \ \mbox{as } n\to \infty,
  \end{eqnarray*}
which implies
\begin{eqnarray}\label{3.2.3}
      \left|
      \left\langle
      \bar{I}_{\lambda}'(u_{n},v_{n}),\left(\frac{1}{\theta_{1}}u_{n},\frac{1}{\theta_{2}}v_{n}\right)
      \right\rangle
      \right|
\leq
       \frac{1}{\min\{\theta_{1},\theta_{2}\}}
       \|\bar{I}_{\lambda}'(u_{n},v_{n})\|_{W^{*}}
        \|(u_{n},v_{n})\|
 \to 0 \ \ \mbox{as } n\to \infty.
 \end{eqnarray}
Then by (\ref{3.2.2}), (\ref{3.2.3}) and $(H_3)'$, for $n$ large enough and any $\lambda>0$, we have
\begin{eqnarray}\label{b13}
         c_{\lambda} +1
 &\geq & \bar{I}_{\lambda}(u_{n},v_{n})
         -\left\langle \bar{I}_{\lambda}'(u_{n},v_{n}),\left(\frac{1}{\theta_{1}}u_{n},\frac{1}{\theta_{2}}v_{n}\right) \right \rangle
          \nonumber \\
 & = &   \left(\frac{1}{p_{1}}-\frac{1}{\theta_{1}}\right)\|u_{n}\|_{m_{1},p_{1}}^{p_{1}}+\left(\frac{1}{p_{2}}-\frac{1}{\theta_{2}}\right)\|v_{n}\|_{m_{1},p_{2}}^{p_{2}}\nonumber \\
 &  &  +\lambda \int_V
 \left(\frac{1}{\theta_{1}}u_{n}\bar{F}_{u}(x,u_{n},v_{n})+\frac{1}{\theta_{2}}v_{n}\bar{F}_{v}(x,u_{n},v_{n}) -\bar{F}(x,u_{n},v_{n})\right)d\mu
           \nonumber \\
 &\geq &  \left(\frac{1}{p_{1}}-\frac{1}{\theta_{1}}\right)\|u_{n}\|_{m_{1},p_{1}}^{p_{1}}+\left(\frac{1}{p_{2}}-\frac{1}{\theta_{2}}\right)\|v_{n}\|_{m_{1},p_{2}}^{p_{2}},
 \end{eqnarray}
which implies that $\{u_{n}\}$ and $\{v_{n}\}$ are bounded in $W^{m_{1},p_{1}}(V)$ and $W^{m_{2},p_{2}}(V)$, respectively.
Due to $W^{m_{1},p_{1}}(V)$ and $W^{m_{2},p_{2}}(V)$ are both pre-compact, up to a subsequence of $\{(u_{n},v_{n})\}$, still denoted by $\{(u_{n},v_{n})\}$, there exists $(u_{\lambda} ,v_{\lambda})\in W$ such that $u_{n}\rightarrow u_\lambda$ in $W^{m_{1},p_{1}}(V)$ and $v_{n}\rightarrow v_\lambda$ in $W^{m_{2},p_{2}}(V)$, respectively.
\qed

\vskip2mm
\par
Lemma \ref{lemma3.4} shows that $(u_{n},v_{n})\rightarrow (u_{\lambda} ,v_{\lambda})$ in $W$, which together with (\ref{3.2.2}) and the continuity of
$\bar{I}_{\lambda}$ and $\bar{I}_{\lambda}'$ imply that $\bar{I}_{\lambda}(u_{\lambda} ,v_{\lambda})=  c_{\lambda}$ and $\bar{I}_{\lambda}'(u_{\lambda} ,v_{\lambda})= 0$. Moreover, by $\bar{I}_{\lambda}(0,0)=  0$ and $c_{\lambda} >0$, we have $(u_{\lambda} ,v_{\lambda})\neq (0,0)$. Thus $(u_{\lambda} ,v_{\lambda})$ is a nontrivial critical point of $\bar{I}_{\lambda}$ with the critical value $c_{\lambda} $. Next, we will show that $(u_{\lambda} ,v_{\lambda})$ precisely is the nontrivial weak solution of system (\ref{aa1}).
\par
\noindent
\begin{lemma}\label{lemma3.5}
Assume that $(H_3)'$ and $(H_5)'$ hold. Then there exist $C_{1,*}>0$ and $C_{2,*}>0$ independent of $\lambda>0$ such that
\begin{eqnarray*}
&  &       \|u_{\lambda}\|_{\infty}
           \leq \left(\dfrac{p_{1}\theta_{1}\max\{C_{1,*},C_{2,*}\}}{(\theta_{1}-p_{1})\mu_{\infty} h_{1,\infty}}\right)^{\frac{1}{p_{1}}}
           \left( \lambda^{-\frac{p_{1}}{q_{1}-p_{1}}}+\lambda^{-\frac{p_{2}}{q_{2}-p_{2}}} \right)^{\frac{1}{p_{1}}},\\
&  &       \|v_{\lambda}\|_{\infty}
           \leq \left( \dfrac{p_{2}\theta_{2}\max\{C_{1,*},C_{2,*}\}}{(\theta_{2}-p_{2})\mu_{\infty} h_{2,\infty}}\right)^{\frac{1}{p_{2}}}
           \left( \lambda^{-\frac{p_{1}}{q_{1}-p_{1}}}+\lambda^{-\frac{p_{2}}{q_{2}-p_{2}}} \right)^{\frac{1}{p_{2}}},
 \end{eqnarray*}
where $C_{1,*}$ and $C_{2,*}$ are defined in (\ref{A6}) and (\ref{A7}), respectively.
\end{lemma}
\vskip2mm
\noindent
{\bf Proof.} Let $s\in [0 , 1]$. Then $\|s(u_{0},v_{0})\|_{\infty}\le \|(u_{0},v_{0})\|_\infty \le  \delta$, where $(u_{0},v_{0})$ is given in  Lemma \ref{lemma3.3}. For all $\lambda>0$, by $(H_5)'$ and H\"older's inequality, we have
\begin{eqnarray}\label{b14}
          \bar{I}_{\lambda}(su_{0},sv_{0})
& = &
         \frac{1}{p_{1}}\|u_{0}\|_{m_{1},p_{1}}^{p_{1}}s^{p_{1}}
          +\frac{1}{p_{2}}\|v_{0}\|_{m_{2},p_{2}}^{p_{2}}s^{p_{2}}
          -\lambda \int_V \bar{F}(x,su_{0},sv_{0})d\mu
\nonumber \\
&\leq &
         \frac{1}{p_{1}}\|u_{0}\|_{m_{1},p_{1}}^{p_{1}}s^{p_{1}}
         +\frac{1}{p_{2}}\|v_{0}\|_{m_{2},p_{2}}^{p_{2}}s^{p_{2}}
         -\lambda M_{1}s^{q_{1}}\int_{V}|u_{0}|^{q_{1}}d\mu
\nonumber \\
&  &
        -\lambda M_{2}s^{q_{2}}\int_{V}|v_{0}|^{q_{2}}d\mu
 \nonumber \\
&\leq &
         \frac{1}{p_{1}}\|u_{0}\|_{m_{1},p_{1}}^{p_{1}}s^{p_{1}}
         -\lambda M_{1}\left(\sum_{x\in V} \mu(x)\right)
         ^{1-\frac{q_{1}}{p_{1}}}
         \|u_{0}\|_{L^{p_{1}}(V)}^{q_{1}}s^{q_{1}}
 \nonumber \\
 &  &
         +\frac{1}{p_{2}}\|v_{0}\|_{m_{2},p_{2}}^{p_{2}}s^{p_{2}}
         -\lambda M_{2}\left(\sum_{x\in V} \mu(x)\right)
         ^{1-\frac{q_{2}}{p_{2}}}
         \|v_{0}\|_{L^{p_{2}}(V)}^{q_{2}}s^{q_{2}}.
  \end{eqnarray}
Define $g_{i}:[0,+\infty)\to \R$, $i=1,2$, by
\begin{eqnarray*}
 &  &      g_{1}(s)= \frac{1}{p_{1}}\|u_{0}\|_{m_{1},p_{1}}^{p_{1}}s^{p_{1}}-\lambda M_{1}\left(\sum_{x\in V} \mu(x)\right) ^{1-\frac{q_{1}}{p_{1}}}\|u_{0}\|_{L^{p_{1}}(V)}^{q_{1}}s^{q_{1}},
 \\
 &  &      g_{2}(s)= \frac{1}{p_{2}}\|v_{0}\|_{m_{2},p_{2}}^{p_{2}}s^{p_{2}}-\lambda M_{2}\left(\sum_{x\in V} \mu(x)\right) ^{1-\frac{q_{2}}{p_{2}}}\|v_{0}\|_{L^{p_{2}}(V)}^{q_{2}}s^{q_{2}}
 .
 \end{eqnarray*}
For any $\lambda>0$, let
$$
g_{1}^{\prime}(s)= s^{p_{1}-1}\left(\|u_{0}\|_{m_{1},p_{1}}^{p_{1}}-\lambda q_{1}M_{1}\left(\sum_{x\in V} \mu(x)\right) ^{1-\frac{q_{1}}{p_{1}}}\|u_{0}\|_{L^{p_{1}}(V)}^{q_{1}}s^{q_{1}-p_{1}} \right)=0.
$$
Then we have
$ s_{1}^{\lambda} = \left(\sum\limits_{x\in V} \mu(x)\right) ^{\frac{1}{p_{1}}} \left(  \dfrac{\|u_{0}\|_{m_{1},p_{1}}^{p_{1}} }{\lambda q_{1}M_{1}\|u_{0}\|_{L^{p_{1}}(V)}^{q_{1}}} \right)^{\frac{1}{q_{1}-p_{1}}} >0$. It follows from Remark \ref{remark2.1}  that
\begin{eqnarray}\label{b15}
\max_{s\geq0}g_{1}(s)
& = & g_{1}(s_{1}^{\lambda})
\nonumber \\
& = & \frac{(q_{1}-p_{1})\sum\limits_{x\in V} \mu(x)}{ p_{1} }
\left(\frac{\|u_{0}\|_{m_{1},p_{1}}}{\|u_{0}\|_{L^{p_{1}}(V)}} \right)^{\frac{p_{1}q_{1}}{q_{1}-p_{1}}}( q_{1}^{q_{1}} M_{1} ^{p_{1}})^{\frac{1}{p_{1}-q_{1}}}
 \lambda^{-\frac{p_{1}}{q_{1}-p_{1}}}\nonumber \\
&\leq & \frac{(q_{1}-p_{1}) \sum\limits_{x\in V}\mu(x) }{p_{1}}
      \left(\frac{ (\|\nabla^{m_{1}} u_{0}\|_{L^{p_{1}}(V)}^{p_{1}}+ \|u_{0}\|_{L^{p_{1}}(V)}^{p_{1}})^{\frac{1}{p_{1}}}} { \|u_{0}\|_{L^{p_{1}}(V)}       }\right)^{\frac{p_{1}q_{1}}{q_{1}-p_{1}}} \nonumber \\
&  &  \left(\frac{(\max\{1,h_{1}^{\infty}\})^{q_{1}}}{q_{1} ^{q_{1}} M_{1} ^{p_{1}}}\right)^{\frac{1}{q_{1}-p_{1}}} \lambda^{-\frac{p_{1}}{q_{1}-p_{1}}}
 \nonumber \\
& = & C_{1,*} \lambda^{-\frac{p_{1}}{q_{1}-p_{1}}}.
 \end{eqnarray}
Similarly, we have
\begin{eqnarray}\label{b16}
\max_{s\geq0}g_{2}(s)
 =  g_{2}(s_{2}^{\lambda})
\leq   C_{2,*} \lambda^{-\frac{p_{2}}{q_{2}-p_{2}}},
\end{eqnarray}
where
$ s_{2}^{\lambda} =  \left(\sum\limits_{x\in V} \mu(x)\right) ^{\frac{1}{p_{2}}}\left(  \dfrac{\|v_{0}\|_{m_{2},p_{2}}^{p_{2}} }{\lambda q_{2}M_{2}\|v_{0}\|_{L^{p_{2}}(V)}^{q_{2}}} \right)^{\frac{1}{q_{2}-p_{2}}} >0$.
By the definition of $c_{\lambda}$, (\ref{b14}), (\ref{b15}) and (\ref{b16}), we have
\begin{eqnarray}\label{b17}
      c_{\lambda}
\leq
     \max_{s \in[0,1]}\bar{I}_{\lambda}(su_{0},sv_{0})
\leq
      \max_{s\geq0}g_{1}(s)+\max_{s\geq0}g_{2}(s)
\leq
      \max\{C_{1,*},C_{2,*}\}
      \left(
      \lambda^{-\frac{p_{1}}{q_{1}-p_{1}}}
      +\lambda^{-\frac{p_{2}}{q_{2}-p_{2}}}
      \right).
\end{eqnarray}
Note that
$\left\langle \bar{I}_{\lambda}'(u_{\lambda} ,v_{\lambda}),(\frac{1}{\theta_{1}}u_{\lambda},\frac{1}{\theta_{2}}v_{\lambda}) \right\rangle =0$.
 Similar to the argument in (\ref{b13}), we have
\begin{equation*}
  \|u_{\lambda}\|_{m_{1},p_{1}}
  \leq \left(\frac{p_{1}\theta_{1}}{\theta_{1}-p_{1}}\right)^{\frac{1}{p_{1}}} c_{\lambda}^{\frac{1}{p_{1}}},\,\,\
  \|v_{\lambda}\|_{m_{2},p_{2}}
  \leq \left(\frac{p_{2}\theta_{2}}{\theta_{2}-p_{2}}\right)^{\frac{1}{p_{2}}} c_{\lambda}^{\frac{1}{p_{2}}},
\end{equation*}
which combining with (\ref{b17}) connote that
\begin{eqnarray}\label{b18}
  \|u_{\lambda}\|_{m_{1},p_{1}}
           \leq \left(\dfrac{p_{1}\theta_{1}\max\{C_{1,*},C_{2,*}\}}{\theta_{1}-p_{1}}\right)^{\frac{1}{p_{1}}}
           \left( \lambda^{-\frac{p_{1}}{q_{1}-p_{1}}}+\lambda^{-\frac{p_{2}}{q_{2}-p_{2}}} \right)^{\frac{1}{p_{1}}}
\end{eqnarray}
and
\begin{eqnarray}\label{b19}
 \|v_{\lambda}\|_{m_{2},p_{2}}
           \leq \left( \dfrac{p_{2}\theta_{2}\max\{C_{1,*},C_{2,*}\}}{\theta_{2}-p_{2}}\right)^{\frac{1}{p_{2}}}
           \left( \lambda^{-\frac{p_{1}}{q_{1}-p_{1}}}+\lambda^{-\frac{p_{2}}{q_{2}-p_{2}}} \right)^{\frac{1}{p_{2}}}.
\end{eqnarray}
By Lemma \ref{lemma2.1} , we further get the following results
\begin{eqnarray}\label{b20}
\|u_{\lambda}\|_{\infty}
           \leq \left(\dfrac{p_{1}\theta_{1}\max\{C_{1,*},C_{2,*}\}}{(\theta_{1}-p_{1})\mu_{\infty} h_{1,\infty}}\right)^{\frac{1}{p_{1}}}
           \left( \lambda^{-\frac{p_{1}}{q_{1}-p_{1}}}+\lambda^{-\frac{p_{2}}{q_{2}-p_{2}}} \right)^{\frac{1}{p_{1}}}
\end{eqnarray}
and
\begin{eqnarray}\label{b21}
\|v_{\lambda}\|_{\infty}
           \leq \left( \dfrac{p_{2}\theta_{2}\max\{C_{1,*},C_{2,*}\}}{(\theta_{2}-p_{2})\mu_{\infty} h_{2,\infty}}\right)^{\frac{1}{p_{2}}}
           \left( \lambda^{-\frac{p_{1}}{q_{1}-p_{1}}}+\lambda^{-\frac{p_{2}}{q_{2}-p_{2}}} \right)^{\frac{1}{p_{2}}}.
\end{eqnarray}
The proof is completed.\qed
\vskip2mm
\noindent
{\bf Proof of Theorem \ref{theorem1.1}}\ \ It follows from the fact $q_{i}>p_{i}>1$, $i=1,2$ and Lemma \ref{lemma2.3} that
$$
     \lambda^{-\frac{p_{1}}{q_{1}-p_{1}}}+\lambda^{-\frac{p_{2}}{q_{2}-p_{2}}}
\leq
     2
     \left(
     1+ \lambda^{-\max\left\{\frac{p_{1}}{q_{1}-p_{1}} ,\frac{p_{2}}{q_{2}-p_{2}} \right\}}
     \right)
$$
which together with Lemma \ref{lemma3.5} implies that $\|(u_{\lambda} ,v_{\lambda})\|_{\infty}\leq \|u_{\lambda}\|_{\infty}+\|v_{\lambda}\|_{\infty}\leq \frac{\delta}{4}+\frac{\delta}{4} = \frac{\delta}{2}$ for each $\lambda>\Lambda_5$, where $\Lambda_5$ is defined in (\ref{A5}). Therefore $(u_{\lambda} ,v_{\lambda})\in W$ is a nontrivial weak solution of  system (\ref{aa1}). Furthermore, (\ref{b18})-(\ref{b21}) imply that
\begin{eqnarray*}
  \|(u_{\lambda},v_{\lambda})\|
    & \leq &      \left(\dfrac{p_{1}\theta_{1}\max\{C_{1,*},C_{2,*}\}}{\theta_{1}-p_{1}}\right)^{\frac{1}{p_{1}}}
           \left( \lambda^{-\frac{p_{1}}{q_{1}-p_{1}}}+\lambda^{-\frac{p_{2}}{q_{2}-p_{2}}} \right)^{\frac{1}{p_{1}}}\\
   &  &  + \left( \dfrac{p_{2}\theta_{2}\max\{C_{1,*},C_{2,*}\}}{\theta_{2}-p_{2}}\right)^{\frac{1}{p_{2}}}
           \left( \lambda^{-\frac{p_{1}}{q_{1}-p_{1}}}+\lambda^{-\frac{p_{2}}{q_{2}-p_{2}}} \right)^{\frac{1}{p_{2}}}
  \end{eqnarray*}
and
 \begin{eqnarray*}
      \|(u_{\lambda},v_{\lambda})\|_{\infty}
    &  \leq &  \left(\dfrac{p_{1}\theta_{1}\max\{C_{1,*},C_{2,*}\}}{(\theta_{1}-p_{1})\mu_{\infty} h_{1,\infty}}\right)^{\frac{1}{p_{1}}}
           \left( \lambda^{-\frac{p_{1}}{q_{1}-p_{1}}}+\lambda^{-\frac{p_{2}}{q_{2}-p_{2}}} \right)^{\frac{1}{p_{1}}}\\
  &  &   +  \left( \dfrac{p_{2}\theta_{2}\max\{C_{1,*},C_{2,*}\}}{(\theta_{2}-p_{2})\mu_{\infty} h_{2,\infty}}\right)^{\frac{1}{p_{2}}}
           \left( \lambda^{-\frac{p_{1}}{q_{1}-p_{1}}}+\lambda^{-\frac{p_{2}}{q_{2}-p_{2}}} \right)^{\frac{1}{p_{2}}}.
 \end{eqnarray*}
Thus, we have
$\lim\limits_{\lambda\to \infty}\|(u_{\lambda},v_{\lambda})\|=0=\lim\limits_{\lambda\to \infty} \|(u_{\lambda},v_{\lambda})\|_\infty. $
\qed
\vskip2mm
 {\section{Proofs for the sub-linear case}}
  \setcounter{equation}{0}
\par
Define $\rho(t,s)\in C^1(\R ^{2},[0,1])$ as a cut-off function satisfying
 \begin{eqnarray*}
 \rho(t,s)= \begin{cases}
           1, \;\;\;  \text { if }\;\;|(t,s)| \leq \delta/2, \\
           0, \;\;\;  \text { if }\;\;|(t,s)|    >   \delta.
          \end{cases}
   \end{eqnarray*}
Define $\tilde{F}:V\times\R\times \R\to\R$ by
   $$
   \tilde{F}(x,t,s)=\rho(t,s)F(x,t,s)+(1-\rho(t,s))(K_{1}|t|^{q_{1}}+K_{2}|s|^{q_{2}}) .
   $$
By $(H_0)$, $(R_2)$-$(R_4)$ and the definition of $\rho(t,s)$, we could obtain the following lemma:
\par
\noindent
\begin{lemma}\label{lemma4.1}
Assume that $(H_0)$, $(R_2)$-$(R_4)$ hold. Then
\begin{itemize}
\item[$(R_0)'$]$\tilde{F}(x,t,s)$ is continuously differentiable in $\R\times \R$ for all $x \in V$;
\item[$(R_2)'$]there exist $K_{7}>0$ and $K_{8}>0$ such that
$$  |\tilde{F}_{t}(x,t,s)|\leq K_{7}(|t|^{q_{1}-1}+|t|^{k_{1}-1})+K_{8}|s|^{\frac{k_{2}(k_{1}-1)}{k_{1}}},$$
$$  |\tilde{F}_{s}(x,t,s)|\leq K_{7}|t|^{\frac{k_{1}(k_{2}-1)}{k_{2}}}+K_{8}(|s|^{q_{2}-1}+|s|^{k_{2}-1})$$
for all $x\in V$ and $(t,s)\in \R ^{2}$;
\item[$(R_3)'$]$\tilde{F}(x,t,s)=\tilde{F}(x,-t,-s)$ for all $x\in V$ and $(t,s)\in \R^{2}$;
\item[$(R_4)'$]
$$K_{1}|t|^{q_{1}}+K_{2}|s|^{q_{2}} \leq \tilde{F}(x,t,s)\leq \max\{K_{1},K_{5}\}(|t|^{q_{1}}+|t|^{k_{1}})+\max\{K_{2},K_{6}\}(|s|^{q_{2}}+|s|^{k_{2}})$$
 for all $x\in V$ and $(t,s)\in \R ^{2}$.
\end{itemize}
\end{lemma}

\par
Consider the modified system of (\ref{aa1}) given by
\begin{equation}\label{mod2}
 \begin{cases}
\pounds_{m_{1},p_{1}}u+h_{1}(x)|u|^{p_{1}-2}u=\lambda \tilde{F}_{u}(x,u,v),  \;\;\; x \in V,\\
\pounds_{m_{2},p_{2}}v+h_{2}(x)|v|^{p_{2}-2}v=\lambda \tilde{F}_{v}(x,u,v),  \;\;\; x \in V.
   \end{cases}
 \end{equation}
 Define the corresponding functional $\tilde{I}_{\lambda}: W\rightarrow \R$ by
 \begin{eqnarray*}
        \tilde{I}_{\lambda}(u,v)
  &  : =  & \frac{1}{p_{1}}\int_{V}(|\nabla^{m_{1}} u|^{p_{1}}
     +h_{1}(x)|u|^{p_{1}})d\mu + \frac{1}{p_{2}}\int_{V}(|\nabla^{m_{2}}v|^{p_{2}}
     +h_{2}(x)|v|^{p_{2}})d\mu \\
    & &  -\lambda \int_V \tilde{F}(x,u,v)d\mu
  \end{eqnarray*}
for all $ (u, v)\in W $. By $(R_0)'$-$(R_4)'$ and the continuous embeddings
 \begin{equation*}
   W^{m_{i},p_{i}}(V)\hookrightarrow L^{q_{i}}(V),\,\,\,\ W^{m_{i},p_{i}}(V)\hookrightarrow L^{k_{i}}(V),\,\,i=1,2,
   \end{equation*}
standard arguments show that $\tilde{I}_{\lambda}$ is well defined, even and of class $C^{1}$ on $W$, and
 \begin{eqnarray*}
          \langle \tilde{I}_{\lambda}'(u,v),(\phi_{1},\phi_{2})\rangle
  &  =  & \int_V\left( \pounds_{m_{1},p_{1}} u \phi_{1}+h_{1}(x)|u|^{p_{1}-2}u\phi_{1}\right)d\mu\\
  &    & + \int_V\left( \pounds_{m_{2},p_{2}} v\phi_{2}+h_{2}(x)|v|^{p_{2}-2}v \phi_{2}\right)d\mu\\
   &    &-\lambda \int_V \tilde{F}_{u}(x,u,v)\phi_{1}d\mu- \lambda \int_V \tilde{F}_{v}(x,u,v) \phi_{2}d\mu
   \end{eqnarray*}
   for all $ (u,v),(\phi_{1},\phi_{2})\in W$ (for example, see \cite{Yang-Zhang-2024}). Hence
   \begin{eqnarray*}
          \langle \tilde{I}_{\lambda}'(u,v),(u,v)\rangle
    =   \|u\|_{m_{1},p_{1}}^{p_{1}}+\|v\|_{m_{2},p_{2}}^{p_{2}}-\lambda \int_V \tilde{F}_{u}(x,u,v)ud\mu- \lambda \int_V \tilde{F}_{v}(x,u,v)vd\mu
   \end{eqnarray*}
 for all $ (u,v)\in W$. Note that the critical points of $\tilde{I}_{\lambda}$ with $ L^{\infty}$-norm less than or equal to $ \delta/2$ are also critical points of $I_{\lambda}$. Then the weak solutions of system (\ref{mod2}) with $ L^{\infty}$-norm less than or equal to $ \delta/2$ are also weak solutions of problem (\ref{aa1}).
\vskip2mm
\noindent
{\bf Proof of Theorem \ref{theorem1.2}}\ \ Obviously, $\tilde{I}_{\lambda}\in C^{1}(W,\R)$ is even and $\tilde{I}_{\lambda}(0,0)=0$. For any fixed $\lambda>0$, using $(R_4)'$ we get
\begin{eqnarray*}
\tilde{I}_{\lambda}(u,v)
  & \geq &   \frac{1}{p_{1}}\|u\|_{m_{1},p_{1}}^{p_{1}}-\lambda \max\{K_{1},K_{5}\}\int_V (|u|^{q_{1}}+|u|^{k_{1}})d\mu \\
  & &  + \frac{1}{p_{2}}\|v\|_{m_{2},p_{2}}^{p_{2}}-\lambda \max\{K_{2},K_{6}\}\int_V (|v|^{q_{2}}+|v|^{k_{2}})d\mu,
\end{eqnarray*}
which estimate by Lemma \ref{lemma2.2} leads to
\begin{eqnarray*}
        \tilde{I}_{\lambda}(u,v)
& \geq &
       \frac{1}{p_{1}}\|u\|_{m_{1},p_{1}}^{p_{1}}
       -\lambda \max\{K_{1},K_{5}\}\sum_{x\in V} \mu(x)
      \left(
        \frac{1}
        {\mu_{\infty}^{\frac{q_{1}}{p_{1}}} h_{1,\infty}^{\frac{q_{1}}{p_{1}}}}
         \|u\|_{m_{1},p_{1}} ^{q_{1}}
        +
        \frac{1}{\mu_{\infty}^{\frac{k_{1}}{p_{1}}} h_{1,\infty}^{\frac{k_{1}}{p_{1}}}}
        \|u\|_{m_{1},p_{1}} ^{k_{1}}
        \right)
\\
& &
        +\frac{1}{p_{2}}\|v\|_{m_{2},p_{2}}^{p_{2}}
        -\lambda \max\{K_{2},K_{6}\}\sum_{x\in V} \mu(x)
       \left(
       \frac{1}{\mu_{\infty}^{\frac{q_{2}}{p_{2}}}h_{2,\infty}^{\frac{q_{2}}{p_{2}}}}
        \|v\|_{m_{2},p_{2}} ^{q_{2}}
         +\frac{1}{\mu_{\infty}^{\frac{k_{2}}{p_{2}}} h_{2,\infty}^{\frac{k_{2}}{p_{2}}}}
         \|v\|_{m_{2},p_{2}} ^{k_{2}}
         \right),
\end{eqnarray*}
which together with the fact $p_{i}>q_{i}>k_{i}>1$, $i=1,2$ implies that $\tilde{I}_{\lambda}(u,v)\rightarrow +\infty$ as $\|(u,v)\|\rightarrow \infty$.
That is $\tilde{I}_{\lambda}$ is coercive and then is bounded form below on $W$.
Suppose that $\{(u_{n},v_{n})\}$ is a Palais-Smale sequence of $\tilde{I}_{\lambda}$ on $W$, that is,
\begin{equation*}
   \{\tilde{I}_{\lambda}(u_{n},v_{n})\}\ \ \mbox{is bounded },
\ \
   \tilde{I}_{\lambda}'(u_{n},v_{n}) \to 0\ \ \mbox{as } n\to \infty.
   \end{equation*}
Then the coercivity of $\tilde{I}_{\lambda}$ implies that $\{(u_{n},v_{n})\}$ is bounded on $W$.
Due to $W^{m_{1},p_{1}}(V)$ and $W^{m_{2},p_{2}}(V)$ are both pre-compact, up to a subsequence of $\{(u_{n},v_{n})\}$, still denoted by $\{(u_{n},v_{n})\}$, there exists $(\tilde{u}_{\lambda} ,\tilde{v}_{\lambda})\in W$ such that $u_{n}\rightarrow \tilde{u}_\lambda$ in $W^{m_{1},p_{1}}(V)$ and $v_{n}\rightarrow \tilde{v}_\lambda$ in $W^{m_{2},p_{2}}(V)$, respectively.
Hence the functional $\tilde{I}_{\lambda}$ satisfies the Palais-Smale condition for all $\lambda>0$.
Let $Z^{k}$ be a $k$-dimensional subspace of $W$ for any $k\in\mathbb{N}$.
Then all norms are equivalent in $Z^{k}$.
Hence, there exist positive constants $c_{1}$ and $c_{2}$ such that
\begin{eqnarray}\label{b12}
       \|u\|_{L^{q_{1}}(V)}^{q_{1}}
\geq
       c_{1}\|u\|_{m_{1},p_{1}}^{q_{1}},
\ \
       \|v\|_{L^{q_{2}}(V)}^{q_{2}}
\geq
       c_{2}\|v\|_{m_{2},p_{2}}^{q_{2}}
\end{eqnarray}
for all $(u,v)\in Z^{k}$. For any given $0<\rho_{k}<\min\{1,\rho_{\lambda}\}$,
\begin{eqnarray*}
  \rho_{\lambda}
=
    \left(
    \frac{p_{1}p_{2}}{p_{1}+p_{2}}
    \lambda \min\{K_{1}c_{1},K_{2}c_{2}\}
      2^{1-\max\{q_{1},q_{2}\}}
      \right)
      ^{\frac{1}{\min\{p_{1},p_{2}\}-\max\{q_{1},q_{2}\}}},
\end{eqnarray*}
define
\begin{eqnarray*}
S_{\rho_{k}}= \{(u,v)\in W|\|(u,v)\|=\rho_{k} \}.
\end{eqnarray*}
By $(R_4)'$, (\ref{b12}) and the fact $\min\{p_{1},p_{2}\}>\max\{q_{1},q_{2}\}$,
we have
\begin{eqnarray*}
        \tilde{I}_{\lambda}(u,v)
 & = &
       \frac{1}{p_{1}}\|u\|_{m_{1},p_{1}}^{p_{1}}
       +\frac{1}{p_{2}}\|v\|_{m_{2},p_{2}}^{p_{2}}
       -\lambda \int_V \tilde{F}(x,u,v)d\mu
\\
 & \leq &
      \frac{p_{1}+p_{2}}{p_{1}p_{2}}
      \|(u,v)\|^{\min\{p_{1},p_{2}\}}
       -\lambda  K_{1}\int_V |u|^{q_{1}}d\mu
       -\lambda  K_{2}\int_V |v|^{q_{2}}d\mu
 \\
 & \leq &
        \frac{p_{1}+p_{2}}{p_{1}p_{2}}
      \|(u,v)\|^{\min\{p_{1},p_{2}\}}
      -\lambda  K_{1}c_{1}\|u\|_{m_{1},p_{1}}^{q_{1}}
      -\lambda  K_{2}c_{2}\|v\|_{m_{2},p_{2}}^{q_{2}}
 \\
 & \leq &
       \frac{p_{1}+p_{2}}{p_{1}p_{2}}
      \|(u,v)\|^{\min\{p_{1},p_{2}\}}
      -\lambda  K_{1}c_{1}\|u\|_{m_{1},p_{1}}^{\max\{q_{1},q_{2}\}}
      -\lambda  K_{2}c_{2}\|v\|_{m_{2},p_{2}}^{\max\{q_{1},q_{2}\}}
 \\
 & \leq &
       \frac{p_{1}+p_{2}}{p_{1}p_{2}}
      \|(u,v)\|^{\min\{p_{1},p_{2}\}}
      -\lambda \min\{K_{1}c_{1},K_{2}c_{2}\}
      (\|u\|_{m_{1},p_{1}}^{\max\{q_{1},q_{2}\}}
      +\|v\|_{m_{2},p_{2}}^{\max\{q_{1},q_{2}\}})
 \\
& \leq &
       \frac{p_{1}+p_{2}}{p_{1}p_{2}}
      \|(u,v)\|^{\min\{p_{1},p_{2}\}}
      -\lambda \min\{K_{1}c_{1},K_{2}c_{2}\}
      2^{1-\max\{q_{1},q_{2}\}}
      \|(u,v)\|^{\max\{q_{1},q_{2}\}}
\\
 & < & 0.
\end{eqnarray*}
for any $(u,v)\in Z^{k} \cap S_{\rho_{k}}$, $\lambda>0$. So $\sup_{Z^{k} \cap S_{\rho_{k}}}\tilde{I}_{\lambda}<0$. Hence, Lemma \ref{lemma2.5} and Remark \ref{remark2.6} show that system (\ref{mod2}) has infinitely many nonzero solutions $\{(u_{k}^{\lambda},v_{k}^{\lambda})\}$ such that $\|(u_{k}^{\lambda},v_{k}^{\lambda})\|\rightarrow0$ as $k\rightarrow\infty$. By Lemma \ref{lemma2.1} we further get
\begin{equation*}
   \|(u_{k}^{\lambda},v_{k}^{\lambda})\|_{\infty}
\leq
    \|u_{k}^{\lambda}\|_{\infty}+\|v_{k}^{\lambda}\|_{\infty}
\leq
    \max\left\{
    \frac{1}{\mu_{\infty}^{\frac{1}{p_{1}}}h_{1,\infty}^{\frac{1}{p_{1}}}},
    \frac{1}{\mu_{\infty}^{\frac{1}{p_{2}}} h_{2,\infty}^{\frac{1}{p_{2}}}}
    \right\}
    \|(u_{k}^{\lambda},v_{k}^{\lambda})\|\rightarrow 0 \ \ \mbox{as } k\to \infty.
   \end{equation*}
So there exists a sufficiently large integer $k_{0}>0$ such that $\|(u_{k}^{\lambda},v_{k}^{\lambda})\|_{\infty}\leq \frac{\delta}{2}$ for all $k\geq k_{0}$, which implies that $\{(u_{k}^{\lambda},v_{k}^{\lambda})\}_{k_{0}}^{\infty}$ is a sequence of weak solutions of the original system (\ref{aa1}) for each fixed $\lambda>0$.
\qed
\vskip2mm
 {\section{Examples}}
  \setcounter{equation}{0}
  \par
 In this section, we present two examples which satisfy the conditions of Theorem \ref{theorem1.1} and Theorem \ref{theorem1.2}, respectively.
\vskip2mm
\par
\noindent
\begin{example}\label{example5.1}
Let $G=(V,E)$ is a finite graph, specifically,
\begin{eqnarray*}
&&
   V=\{x_{1}, x_{2}, x_{3}, x_{4}, x_{5}, x_{6}, x_{7}\},
\\
&&
    E=\{
        \{x_{1}, x_{2}\},
        \{x_{1}, x_{3}\},
        \{x_{1}, x_{4}\},
        \{x_{1}, x_{5}\},
        \{x_{2}, x_{3}\},
        \{x_{2}, x_{6}\},
       \{x_{3}, x_{4}\},
       \{x_{4}, x_{7}\},
       \{x_{5}, x_{6}\},
       \{x_{5}, x_{7}\}
      \},
\\
&&
     \mu(x_{1})=2,
\;     \mu(x_{2})=1,
\;     \mu(x_{3})=\mu(x_{4})=3,
\;     \mu(x_{5})=\mu(x_{6})=5,
\;     \mu(x_{7})=6
\;
\;
\mbox{and}
\;
\;
       deg(x_{1})=6.
\end{eqnarray*}
Then $\int_{V}1d\mu =\Sigma_{x\in V}\mu(x)=25$.
Let $A=\{x\in V|x\sim x_{1}\}$ and $\sharp A=4$, where $\sharp A$ is the number of elements in the set $A$.
Let
\begin{equation*}
e_{1}(x)
=e_{2}(x)
=
 \begin{cases}
1, &\mbox{if} \;\;x=x_{1},\\
0, &\mbox{if} \;\;x\neq x_{1}.
 \end{cases}
 \end{equation*}
Choose
\begin{equation}\label{exa5}
  u_{0}
=
  \frac{\delta \mu_{\infty}^{\frac{1}{2}}h_{1,\infty}^{\frac{1}{2}}}
  {2\sqrt{17}(\max\{1,h_{1}^{\infty}\})^{\frac{1}{2}}}
  e_{1}
\;\;
\mbox{and}
\;\;
   v_{0}
=
  \frac{\delta \mu_{\infty}^{\frac{1}{2}}h_{2,\infty}^{\frac{1}{2}}}
  {2\sqrt{17}(\max\{1,h_{2}^{\infty}\})^{\frac{1}{2}}}
  e_{2}.
   \end{equation}
Let $\delta=1$, $m_{1}=m_{2}=1$, $p_{1}=p_{2}=2$ and $h_{1}(x)=h_{2}(x)=1$ for all $x\in V$.
Consider the following system:
\begin{equation}\label{exa1}
 \begin{cases}
-\Delta u+u=\lambda F_{u}(x,u,v),  \;\;\; x \in V,\\
-\Delta v+v=\lambda F_{v}(x,u,v),  \;\;\; x \in V,
   \end{cases}
 \end{equation}
where
\begin{equation*}
 F(x,t,s)= \sigma(t,s)(|t|^{6}+|s|^{6})+(1-\sigma(t,s))(|t|^{4}+|s|^{4})
   \end{equation*}
for all $(x,t,s)\in V\times\R\times\R$ with $\sigma(t,s)$ defined by
\begin{equation*}
\sigma(t,s)
=
 \begin{cases}
1, &\mbox{if} \;\;|(t,s)|\leq1,\\
\sin\frac{\pi(t^{2}+s^{2}-16)^{2}}{450}, &\mbox{if}\;\; 1<|(t,s)|\leq4,\\
0, &\mbox{if} \;\;|(t,s)|>4.
 \end{cases}
 \end{equation*}
So,
\begin{equation*}
F(x,t,s)
=
 \begin{cases}
|t|^{6}+|s|^{6},&\mbox{if} \;\;|(t,s)|\leq1,\\
(|t|^{6}+|s|^{6})\sin\frac{\pi(t^{2}+s^{2}-16)^{2}}{450}\\
+\left( 1- \sin\frac{\pi(t^{2}+s^{2}-16)^{2}}{450}  \right)(|t|^{4}+|s|^{4}), &\mbox{if}\;\; 1<|(t,s)|\leq4,\\
|t|^{4}+|s|^{4}, &\mbox{if} \;\;|(t,s)|>4.
 \end{cases}
 \end{equation*}
By Theorem \ref{theorem1.1}, we can obtain that system  (\ref{exa1}) has at least a nontrivial weak solution $(u_{\lambda},v_{\lambda})$  for each $\lambda>\Lambda_4\approx 89523333.3$  and
$\lim_{\lambda\to \infty} \|(u_{\lambda},v_{\lambda})\|=0=\lim_{\lambda\to \infty} \|(u_{\lambda},v_{\lambda})\|_\infty$.
\par
In fact, we can verify that $F$ satisfies the conditions in Theorem \ref{theorem1.1}.
Let
\begin{eqnarray*}
    A(t,s)
&:=&
   6|t|^{4}t\sin\frac{\pi(t^{2}+s^{2}-16)^{2}}{450}
+
   4\left( 1- \sin\frac{\pi(t^{2}+s^{2}-16)^{2}}{450}  \right)b(x)|t|^{2}t
\\
& &
+
   (|t|^{6}+|s|^{6})\frac{2\pi t(t^{2}+s^{2}-16)}{225}\cos\frac{\pi(t^{2}+s^{2}-16)^{2}}{450}
\\
& &
-
   (|t|^{4}+|s|^{4})\frac{2\pi t(t^{2}+s^{2}-16)}{225}\cos\frac{\pi(t^{2}+s^{2}-16)^{2}}{450},
 \end{eqnarray*}
and
\begin{eqnarray*}
    B(t,s)
&:=&
   6|s|^{4}s\sin\frac{\pi(t^{2}+s^{2}-16)^{2}}{450}
+
   4\left( 1- \sin\frac{\pi(t^{2}+s^{2}-16)^{2}}{450} \right)|s|^{2}s
\\
& &
+
   (|t|^{6}+|s|^{6})\frac{2\pi s(t^{2}+s^{2}-16)}{225}\cos\frac{\pi(t^{2}+s^{2}-16)^{2}}{450}
\\
& &
-
   (|t|^{4}+|s|^{4})\frac{2\pi s(t^{2}+s^{2}-16)}{225}\cos\frac{\pi(t^{2}+s^{2}-16)^{2}}{450}.
 \end{eqnarray*}
It is easy to see that $F$ satisfies $(H_0)$ and
\begin{equation*}
F_{t}(x,t,s)
=
 \begin{cases}
6|t|^{4}t,&\mbox{if} \;\;|(t,s)|\leq1,\\
A(t,s), &\mbox{if}\;\; 1<|(t,s)|\leq4,\\
4|t|^{2}t, &\mbox{if} \;\;|(t,s)|>4,
 \end{cases}
 \end{equation*}
\begin{equation*}
F_{s}(x,t,s)
=
 \begin{cases}
6|s|^{4}s,&\mbox{if} \;\;|(t,s)|\leq1,\\
B(t,s), &\mbox{if}\;\; 1<|(t,s)|\leq4,\\
4|s|^{2}s, &\mbox{if} \;\;|(t,s)|>4.
 \end{cases}
 \end{equation*}
Hence, we get
\begin{eqnarray*}
    F(x,t,s)
=  |t|^{6}+|s|^{6}
\geq
  |t|^{7}+|s|^{7},\;\; \mbox{for all} \;\; x\in V \;\; \mbox{and} \;\; |(t,s)|\leq1.
 \end{eqnarray*}
Then $(H_1)$ holds with $q_{1}=q_{2}=7>2=p_{1}=p_{2}$, $M_{1}=M_{2}=1$. We also have
\begin{eqnarray*}
   |F_{t}(x,t,s)|
=6|t|^{5}
\leq
  6(|t|^{4}+|s|^{4}),\;\; \mbox{for all} \;\; x\in V \;\; \mbox{and} \;\; |(t,s)|\leq1.
 \end{eqnarray*}
Similarly,
\begin{eqnarray*}
   |F_{s}(x,t,s)|
\leq
  6(|t|^{4}+|s|^{4}),\;\; \mbox{for all} \;\; x\in V \;\; \mbox{and} \;\; |(t,s)|\leq1.
 \end{eqnarray*}
Then $(H_2)$ holds with $k_{1}=k_{2}=5\in (2,7)$ and $M_{3}=M_{4}=6$.
In particular, taking $\beta_{1}=\beta_{2}=3>2=p_{1}=p_{2}$, we have
\begin{eqnarray*}
   F(x,t,s)
=
     |t|^{6}+|s|^{6}
\leq
     2(|t|^{6}+|s|^{6})
=
   \frac{1}{\beta_{1}}tF_{t}(x,t,s)+ \frac{1}{\beta_{2}}sF_{s}(x,t,s)
 \end{eqnarray*}
for all $x\in \R^{N}$  and $|(t,s)|\leq1$.
Thus, assumption $(H_3)$ is satisfied.
\par
Next, we complete the value of $\lambda_{\ast}$ by the formulas in Theorem \ref{theorem1.1}.
It follows that
\begin{eqnarray}\label{exa3}
      \|\nabla e_{1}\|_{L^{2}(V)}^{2}
&=&
   \sum_{x\in V}|\nabla e_{1}|^{2}(x) \mu(x)
\nonumber\\
&=&
    \sum_{x\in V}
     \left(
     \frac{1}{2\mu(x)} \sum_{y\sim x} \omega_{xy} (e_{1}(y)-e_{1}(x))^{2}
     \right)
     \mu(x)
\nonumber\\
&=&
     \frac{1}{2\mu(x_{1})} \sum_{y\sim x_{1}} \omega_{x_{1}y} \mu(x_{1})
     +
    \sum_{x\sim x_{1}}
     \left(
     \frac{1}{2\mu(x)} \sum_{x_{1}\sim x} \omega_{xx_{1}}
    \right)
     \mu(x)
\nonumber\\
&=&
    \frac{1}{2}deg(x_{1})(1+\sharp A)
\nonumber\\
&=&
    15
=
     \|\nabla e_{2}\|_{L^{2}(V)}^{2}.
 \end{eqnarray}
 Similarly, we have
\begin{eqnarray}\label{exa4}
      \| e_{1}\|_{L^{2}(V)}^{2}
=
     \|e_{2}\|_{L^{2}(V)}^{2}
=
     \mu(x_{1})=2.
 \end{eqnarray}
By Lemma \ref{lemma2.1}, Remark \ref{remark2.1}, (\ref{exa3}), (\ref{exa4}) and (\ref{exa5}), we have
\begin{equation*}
  0<\|u_{0}\|_{\infty}\leq \frac{\delta}{2}
\;\;
\mbox{and}
\;\;
   0<\|v_{0}\|_{\infty}\leq \frac{\delta}{2}.
   \end{equation*}
Moreover,
\begin{eqnarray*}
&&
     \|\nabla u_{0}\|_{L^{2}(V)}^{2}
    = \frac{15\delta^{2} \mu_{\infty}h_{1,\infty}}
       {68\max\{1,h_{1}^{\infty}\}},
\;\;\;\;
     \|\nabla v_{0}\|_{L^{2}(V)}^{2}
    = \frac{15\delta^{2} \mu_{\infty}h_{2,\infty}}
       {68\max\{1,h_{2}^{\infty}\}},
\\
&&
      \|u_{0}\|_{L^{2}(V)}^{2}
   = \frac{\delta^{2} \mu_{\infty}h_{1,\infty}}
       {34\max\{1,h_{1}^{\infty}\}},
\;\;\;\;
      \|v_{0}\|_{L^{2}(V)}^{2}
    = \frac{\delta^{2} \mu_{\infty}h_{2,\infty}}
       {34\max\{1,h_{2}^{\infty}\}}.
   \end{eqnarray*}
Then by
$\theta_{1}=\theta_{2}=3$,
$M_{5}=\frac{12}{5}$,
$M_{6}=6$,
$h_{1,\infty}=h_{2,\infty}=1=h_{1}^{\infty}=h_{2}^{\infty}$,
$\mu_\infty=\min\limits_{x\in V} \mu(x)=1$,
(\ref{A6}) and (\ref{A7}), we have
\begin{eqnarray*}
C_{1,*}=C_{2,*}=\frac{125}{2}\left(\frac{17}{14}\right)^{\frac{7}{5}},
\end{eqnarray*}
and by (\ref{A1})-(\ref{A5}),
 \begin{eqnarray*}
 & &
\Lambda_1=\frac{5}{168},
\;\;\;
 \Lambda_2=\frac{5}{168}(6)^{-3}(34)^{\frac{3}{2}},
\\
 & &
\Lambda_3=\frac{5}{168}(3)^{-3}(34)^{\frac{3}{2}},
\;\;\;
\Lambda_4=\frac{17}{136}(25)^{\frac{5}{2}} (34)^{\frac{7}{2}},
\\
 & &
\Lambda_5=\frac{(2)^{9}(3)^{\frac{5}{2}}(5)^{\frac{15}{2}} (17)^{\frac{7}{2}}}
           {
           \left(
           (7)^{\frac{7}{5}}-3 (5)^{3}(2)^{\frac{18}{5}}(17)^{\frac{7}{5}}
           \right)^{\frac{5}{2}}}.
\end{eqnarray*}
Compared $ \Lambda_1$, $ \Lambda_2$, $ \Lambda_3$, $ \Lambda_4$ and $ \Lambda_5$, it is easy to see $\lambda_{\ast}=\Lambda_4\approx 89523333.3$.
\end{example}
\vskip2mm
\par
\noindent
\begin{remark}\label{remark5.1}
Through the conclusion of Theorem 1.1 and the calculations for the parameter $\lambda$ in Example 5.1, it is easy to see the following interesting phenomenon: the specific lower bound of the parameter on a finite graph is related to the structure of the graph, such as the degree and the measure of any point in $V$. However, in the Euclidean setting, in general, the specific lower bound is only related to the measure of domain (for example, see \cite{Kang-Liu2020}). The fundamental reason for this distinction lies in the difference between the definitions of integrals on the graph and in Euclidean space.
\end{remark}
\noindent
\begin{example}\label{example5.2}
Let $G=(V,E)$ is a finite graph, $\delta=1$, $m_{1}=m_{2}=1$, $p_{1}=2$, $p_{2}=3$ and $h_{1}(x)=h_{2}(x)=1$ for all $x\in V$.
Consider the following system:
\begin{equation}\label{exa6}
 \begin{cases}
-\Delta u+u=\lambda F_{u}(x,u,v),  \;\;\; x \in V,\\
-\Delta_{3} v+|v|v=\lambda F_{v}(x,u,v),  \;\;\; x \in V,
   \end{cases}
 \end{equation}
where
\begin{equation*}
 F(x,t,s)= \frac{3}{4}(|t|^{\frac{4}{3}}+|s|^{\frac{4}{3}})
   \end{equation*}
for all $(x,t,s)\in V\times\R\times\R$.
By Theorem 1.2, we can obtain that system  (\ref{exa6}) has a sequence of weak solutions $\{(u_{k}^{\lambda},v_{k}^{\lambda})\}$ with $\|(u_{k}^{\lambda},v_{k}^{\lambda})\|\rightarrow 0$ as $k\rightarrow\infty$ for every $\lambda>0$.
\par
In fact, we can verify that $F$ satisfies the conditions in Theorem \ref{theorem1.2}.
It is easy to see that $(H_0)$, $(H_4)$ and $(R_3)$ hold.
We also have $(R_1)$ holds with $q_{1}=q_{2}=\frac{5}{3}\in (1,2)$ and  $K_{1}=K_{2}=\frac{3}{4}$.
It follows that
\begin{eqnarray*}
   |F_{t}(x,t,s)|
=|t|^{\frac{1}{3}}
\leq
  2(|t|^{\frac{1}{4}}+|s|^{\frac{1}{4}}),\;\; \mbox{for all} \;\; x\in V \;\; \mbox{and} \;\; |(t,s)|\leq1.
 \end{eqnarray*}
Similarly,
\begin{eqnarray*}
   |F_{s}(x,t,s)|
\leq
  2(|t|^{\frac{1}{4}}+|s|^{\frac{1}{4}}),\;\; \mbox{for all} \;\; x\in V \;\; \mbox{and} \;\; |(t,s)|\leq1.
 \end{eqnarray*}
Then $(R_2)$ holds with $k_{1}=k_{2}=\frac{5}{4}\in (1,\frac{5}{3})$ and $K_{3}=K_{4}=2$.
\end{example}

 \vskip2mm
 \noindent
 {\bf Conflict of interest}\\
 The authors declare that they have no conflict of interest.
 \vskip2mm
 \noindent
 {\bf Funding}\\
This project is supported by Yunnan Fundamental Research Projects (grant No: 202301AT070465) and  Xingdian Talent Support Program for Young Talents of Yunnan Province.

\vskip2mm
 \noindent
 {\bf Authors' contributions}\\
Wanting Qi and Xingyong Zhang contribute the main manuscript equally.

\vskip2mm
 \noindent
 {\bf Data availability statement}\\
 No data is used in this paper.
 \vskip2mm
 \noindent
 {\bf Ethical approval}\\
Not applicable.

\renewcommand\refname{References}
{}
 \end{document}